\documentclass{article}

\usepackage{amsmath}
\usepackage{amssymb}
\usepackage{graphicx}
\usepackage{graphics}
\usepackage{amsthm}
\usepackage{amscd}
\usepackage{color}
\usepackage{amsfonts}
\usepackage{fancyhdr}
\usepackage{galois}
\usepackage{fullpage}
\newtheorem{theorem}{Theorem}[section]
\newtheorem{thm}{Theorem}
\newtheorem{mainthm}{Theorem}
\newtheorem*{theorem*}{Theorem}
\newtheorem{corollary}[theorem]{Corollary}
\newtheorem{proposition}[theorem]{Proposition}
\newtheorem{lemma}[theorem]{Lemma}
\newtheorem{remark}[theorem]{Remark}
\newtheorem{claim}[theorem]{Claim}
\newtheorem*{definition*}{Definition}

\newtheorem*{Question*}{Question}
\newtheorem{Question}{Question}

\theoremstyle{definition}
\newtheorem{definition}[theorem]{Definition}

\numberwithin{equation}{section}

\def\N{\mathbb{N}}

\def\norm #1{\Vert \,#1\, \Vert\,}
\makeatletter

\newcommand{\Rmnum}[1]{\expandafter\@slowromancap\romannumeral #1@}
\makeatother

\def \diff {\operatorname{Diff}}
\def \dim {\operatorname{dim}}
\def \orb {\operatorname{Orb}}
\def \ind {\operatorname{Ind}}
\def \Jac {\operatorname{Jac}}

\def\ud{\mathrm{d}}
\def\um{\mathrm{m}}
\begin{document}

\vspace{-2cm}
%%%%%%%%%%%%%%%%%%%%%%%%%%%%%%%%%%%%%%%%%%%%%%%%%%%%%%%%%%%%%%
%%%%%%%%%%%%%%%%%%%%%%%%%%%%%%%%%%%%%%%%%%%%%%%%%%%%%%%%%%%%%%%%
\title{Ergodic measures with multi-zero Lyapunov exponents inside homoclinic classes}

\author{Xiaodong Wang and Jinhua Zhang}
%%%%%%%%%%%%%%%%%%%%%%%%%%%%%%%%%%%%%%%%%%%%%%%%%%%%%%%%%%%%%%%%%%
%%%%%%%%%%%%%%%%%%%%%%%%%%%%%%%%%%%%%%%%%%%%%%%%%%%%%%%%%%%%%%%%%%%5

\vspace{-2cm}
\maketitle

\begin{abstract}

We prove that for $C^1$ generic diffeomorphisms, if a homoclinic class $H(P)$ contains two hyperbolic periodic orbits of indices $i$ and $i+k$ respectively and $H(P)$ has no domination of index $j$  for any $j\in\{i+1,\cdots,i+k-1\}$, then there exists a non-hyperbolic ergodic measure whose $(i+l)^{th}$ Lyapunov exponent  vanishes  for any $l\in\{1,\cdots, k\}$, and whose support is the whole homoclinic class.

We also prove that for $C^1$ generic diffeomorphisms, if a homoclinic class $H(P)$ has a dominated splitting of the form  $E\oplus F\oplus G$,
such that the center bundle $F$ has no finer dominated splitting, and $H(p)$ contains a hyperbolic periodic orbit $Q_1$ of index $\dim(E)$ and a hyperbolic periodic orbit $Q_2$ whose absolute Jacobian along the bundle $F$ is strictly less than $1$,
then there exists a non-hyperbolic ergodic measure whose Lyapunov exponents along the center bundle $F$ all vanish  and
whose support is the whole homoclinic class.
\end{abstract}

%%%%%%%%%%%%%%%%%%%%%%%%%%%%%%%%%%%%%%%%%%%%%%%%%%%%%%%%%%%%%%%%%%%%%%%%%%%%
%%%%%%%%%%%%%%%%%%%%%%%%%%%%%%%%%%%%%%%%%%%%%%%%%%%%%%%%%%%%%%%%%%%%%%%%%%%%
\section{Introduction}
%%%%%%%%%%%%%%%%%%%%%%%%%%%%%%%%%%%%%%%%%%%%%%%%%%%%%%%%%%%%%%%%%%%%%%%%%%%%%
%%%%%%%%%%%%%%%%%%%%%%%%%%%%%%%%%%%%%%%%%%%%%%%%%%%%%%%%%%%%%%%%%%%%%%%%%%%%%
Since the middle of last century, the dynamics of hyperbolic systems are well understood by dynamicists. Hyperbolic systems have many good properties, for example $\Omega$-stability and existence of Markov partition. However, it was shown by R. Abraham and S. Smale~\cite{AS}  that the hyperbolic systems are not dense among all the differential dynamical systems.
Pesin's theory~\cite{pesin} gives a new notation of hyperbolicity called non-uniform hyperbolicity, which also exhibits asymptotic expansion and contraction rate on the tangent space but may not have uniform bounds for the expansion and contraction time. The example by~\cite{CLR} shows that there exists a non-uniform hyperbolic system exhibiting homoclinic tangencies. Hence, non-uniformly hyperbolic system in general is not hyperbolic. Nevertheless, a series of works by Y. Pesin and A.Katok (for example~\cite{K} and \cite{pesin}) show that many good properties of hyperbolic systems would survive in the non-uniformly hyperbolic setting, for example shadowing property and existence of stable and unstable manifolds. Then, it's natural to ask if the non-uniformly hyperbolic systems are dense among all the differential systems. The first counterexample was given by \cite{KN} in a global setting (some special partially hyperbolic diffeomorphisms), showing that the existence of non-hyperbolic ergodic measures is persistent. Recently, another example is given by \cite{BBD2} in a local setting.
\medskip

Let $M$ be a smooth compact Riemannian manifold of dimension $d$ without boundary. Denote by $\diff^1(M)$ the space of $C^1$ diffeomorphisms of $M$. Consider a diffeomorphism $f\in \diff^1(M)$. By Oseledets's Theorem \cite{O}, for an $f$-invariant ergodic measure $\nu$, there exist $d$ numbers $\chi_{1}(\nu,f)\leq\chi_{2}(\nu,f)\leq\cdots\leq\chi_{d}(\nu,f)$ and a $\nu$-full measure set $\Lambda$ which is invariant under $f$, satisfying that for any $x\in\Lambda$ and any vector $v\in T_x M\setminus\{0\}$, there exists $i\in\{1,2\cdots,d\}$ such that
$$\lim_{k\rightarrow\infty}\frac{1}{k}\log\norm{{\rm Df}^{k}v}=\chi_{i}(\nu,f).$$
The number $\chi_i(\nu,f)$ is called the $i^{th}$ \emph{Lyapunov exponent} of $\nu$. The measure $\nu$ is called \emph{hyperbolic}, if all of its Lyapunov exponents are non-zero. In particular, if $\nu$ is an  atomic measure distributed averagely on a periodic orbit $P=\orb(p)$, then its $i^{th}$ Lyapunov exponent is also called the \emph{$i^{th}$ Lyapunov exponent} of $P$ and is denoted by $\chi_i(p,f)$ or $\chi_i(P,f)$. Assume $E$ is a $\rm{Df}$-invariant subbundle of $T_{\Lambda}M$, then the Lyapunov exponents corresponding to the vectors in $E$ are called the \emph{Lyapunov exponents along $E$}. The number of negative Lyapunov exponents of a hyperbolic measure $\nu$ (or a hyperbolic periodic orbit $P$) is called the \emph{index} of $\nu$ (or $P$), denoted by $\ind(\nu)$ (or $\ind(P)$).
\medskip

The dynamics of a system essentially concentrates on the set of points that have some recurrence properties, the \emph{chain recurrent set} for instance, which splits into disjoint invariant compact sets called \emph{chain recurrence classes}. By~\cite{BC}, for $C^1$-generic diffeomorphisms (i.e. diffeomorphisms in a dense $G_{\delta}$ subset of $\diff^1(M)$), the chain recurrent set coincides with the closure of the set of periodic points and each chain recurrence class containing a periodic orbit $P=\orb(p)$ coincides with its \emph{homoclinic class} $H(P,f)$: the closure of the transverse intersections of the stable and unstable manifolds of $P$.
\medskip

Given an invariant compact set $\Lambda$. We say that $\Lambda$ admits a \emph{$T$-dominated splitting} for a positive integer $T$, if the tangent bundle has a non-trivial ${\rm Df}$-invariant splitting $T_{\Lambda}M=E\oplus F$ such that $$\norm{{\rm Df}^T|_{E(x)}}\norm{{\rm Df}^{-T}|_{F(f^T(x))}}<\frac{1}{2}, \textrm{ for any $x\in\Lambda$}.$$ We say $\Lambda$ admits a \emph{dominated splitting}, if it admits a $T$-dominated splitting for some positive integer $T$. The dimension of the bundle $E$ is called the \emph{index of the dominated splitting}.
\medskip

Recall that a property is called a \emph{generic} property if it is satisfied for a dense $G_{\delta}$ subset of $\diff^1(M)$. There are some previous works to characterize the non-hyperbolicity of homoclinic classes by the existence of non-hyperbolic ergodic measures supported on it, for example~\cite{DG,BDG,CCGWY}. Some method is introduced in~\cite{GIKN} to obtain the ergodicity of  weak-$*$-limit measure of atomic measures supported on periodic orbits, and it is developed in~\cite{DG,BDG}.

\begin{thm}[\cite{DG,BDG}]\label{Thm:D-G-B}
For generic $f\in\diff^1(M)$, consider a hyperbolic periodic orbit $P$ of index $i$. Assume that the homoclinic class $H(P,f)$ contains a hyperbolic periodic orbit $Q$ of index $i-1$, then $H(P,f)$ supports a non-hyperbolic ergodic measure, whose $i^{th}$ Lyapunov exponent vanishes.

 If moreover  $H(P,f)$ admits a dominated splitting $T_{H(P,f)}M=E\oplus F\oplus G$ with $\dim(E)=i-1$ and $\dim(F)=1$, then $H(P,f)$ supports a non-hyperbolic ergodic measure whose Lyapunov exponent along the bundle $F$ vanishes and whose support equals $H(P,f)$.
\end{thm}

Based on the results of~\cite{DG,BDG} and combined to the results of~\cite{bcdg,wang}, a recent work of~\cite{CCGWY} shows that  for $C^1$-generic diffeomorphisms, if a homoclinic class is not hyperbolic, then it supports a non-hyperbolic ergodic measure. Moreover, if the homoclinic class contains periodic orbits of different indices, then one can obtain a non-hyperbolic ergodic measure whose support is the whole homoclinic class.
\medskip

The non-hyperbolic ergodic measures in the discussions above can be only assured to have one vanishing Lyapunov exponent. The example in~\cite{BBD1} shows that there exist iterated function systems (IFS) persistently exhibiting non-hyperbolic ergodic measures with all the Lyapunov exponents vanished. Here, we restate the question posed in \cite{BBD1}:

\begin{Question}\label{Que:multiple 0 exponents}
Does there exist an open set $\mathcal{U}$ of diffeomorphisms such that for any $f\in\mathcal{U}$, there exists an ergodic measure with more than one vanishing Lyapunov exponents~?
\end{Question}

Also one can ask a similar question for homoclinic classes.

\begin{Question}\label{Que:homoclinic class}
Under what kind of assumption, does there exist a non-hyperbolic ergodic measure supported on a homoclinic class with more than one vanishing Lyapunov exponents~?
\end{Question}

Inspired by Theorem~\ref{Thm:D-G-B}, we would like to consider the question that:\\
\textit{If  a  homoclinic class contains periodic points of indices $i$ and $i+k$ respectively, $k>0$, does there exist an ergodic measure supported on the homoclinic class such that all its $(i+1)^{th}$ to $(i+k)^{th}$ Lyapunov exponents vanish~?}

Obviously, it is not true if $H(P,f)$ admits a dominated splitting of index $i+j$ for some $1\leq j\leq k-1$. What happens when there is no such dominated splitting over the class~? We state our first result, which partially answers Question~\ref{Que:homoclinic class}.

\begin{mainthm}\label{Thm: A} For generic diffeomorphism $f\in\diff^1(M)$, consider a hyperbolic periodic orbit $P$. Assume the homoclinic class $H(P,f)$ satisfies the following properties:
\begin{itemize}
\item[--] $H(P,f)$ contains hyperbolic periodic orbits of indices $i$ and $i+k$ respectively, where $i,k>0$;
\item[--] for any integer $1\leq j\leq k-1$, there is no dominated splitting of index $i+j$ over $H(P,f)$.
\end{itemize}
Then there exists an ergodic measure $\nu$ whose support is $H(P,f)$ such that the $(i+j)^{th}$ Lyapunov exponent of $\mu$ vanishes for any $1\leq j\leq k$.
\end{mainthm}

\begin{remark}
Considering the support of the non hyperbolic ergodic measure, in Theorem~\ref{Thm: A}, the case when $k=1$ can be obtained as a combination of Theorem~\ref{Thm:D-G-B} above and Theorem B of~\cite{CCGWY}: if there is a dominated splitting into three bundles, then one can apply Theorem~\ref{Thm:D-G-B}; otherwise Theorem B of~\cite{CCGWY} concludes.
\end{remark}

One has the following direct corollary of Theorem~\ref{Thm: A}, which generalizes the ``moreover'' part of Theorem~\ref{Thm:D-G-B} in the sense that one can obtain non-hyperbolic ergodic measure with more than one vanishing Lyapunov exponents.

\begin{corollary} \label{thm.with dominated splitting}
For generic diffeomorphism $f\in\diff^1(M)$, consider a hyperbolic periodic orbit $P$. Assume that the homoclinic class $H(P,f)$ has a dominated splitting $T_{H(P,f)}M=E\oplus F\oplus G$. Assume, in addition, that the followings are satisfied:
\begin{itemize}
\item[--] $H(P,f)$ contains hyperbolic periodic orbits of indices $\dim(E)$ and $\dim(E\oplus F)$ respectively,
\item[--] the center bundle $F$ has no finer dominated splitting.
\end{itemize}
Then there exists an ergodic measure $\nu$ whose Lyapunov exponents along the bundle $F$ vanish, and whose support is $H(P,f)$.
\end{corollary}

%It can be seen that, in the setting of Corollary~\ref{thm.with dominated splitting}, the case when $\dim(F)=1$
%is the ``moreover'' part of Theorem~\ref{Thm:D-G-B}.

We point out that the assumption of existence of both periodic orbits of indices $\dim(E)$ and $\dim(E\oplus F)$ is important. We can give an example based on the results of~\cite{B,BV}, showing that if there is no periodic orbit of index $\dim(E\oplus F)$ inside the homoclinic class, the conclusion of  Corollary~\ref{thm.with dominated splitting} may not be valid. Actually, in the example, the center bundle $F$ has no finer domination but $F$ is uniformly volume expanding, which forbids to have non-hyperbolic ergodic measures with all zero center Lyapunov exponents
(See the details in Section~\ref{Section:counterexample}). One can also ask the following question, to consider the case when the center bundle $F$ is not volume expanding.

\begin{Question}\label{Q:Jacobian}
In the assumption of \emph{Corollary~\ref{thm.with dominated splitting}}, if we replace the existence of hyperbolic periodic orbit  of index $\dim(E\oplus F)$ by the existence of hyperbolic periodic orbit  whose absolute Jacobian along center bundle $F$ is strictly less than $1$,  does there exist an ergodic measure $\nu$ supported on $H(P,f)$ such that all the Lyapunov exponents of $\nu$ along $F$ vanish~?
\end{Question}

The following theorem gives an affirmative answer to Question \ref{Q:Jacobian}. For a periodic orbit $Q=\orb(q)$, we denote by $\pi(Q)$ (or $\pi(q)$) its period.

\begin{mainthm}\label{Thm: B} For generic diffeomorphism $f\in\diff^1(M)$, consider a hyperbolic periodic orbit $P$. Assume that  the homoclinic class $H(P,f)$ admits a dominated splitting $T_{H(P,f)}M=E\oplus F\oplus G$. Assume, in addition, that we have the following:
\begin{itemize}
\item[--] $H(P,f)$ contains a hyperbolic  periodic orbit of index $\dim(E)$ and a hyperbolic periodic point $q\in H(P,f)$ such that $$|\Jac({\rm Df}^{\pi(q)}|_{F(q)})|<1;$$
\item[--] the center bundle $F$ has no finer dominated splitting.
\end{itemize}
Then there exists an ergodic measure $\nu$ whose support is $H(P,f)$, such that all the Lyapunov exponents of $\nu$ along $F$ vanish.
\end{mainthm}
\begin{remark}
$(1)$ It's clear that \emph{Corollary~\ref{thm.with dominated splitting}} is also implied by Theorem \ref{Thm: B}.

$(2)$ We point out here that, under the assumption of Theorem~\ref{Thm: B}, by applying Theorem 1 of~\cite{bcdg} inductively, one can obtain that $H(P,f)$ contains periodic points with indices equal to $i+k-1$ whose $(i+k)^{th}$ Lyapunov exponent (positive but) arbitrarily close to 0.
\end{remark}

Let's explain a little bit about the relation between Theorem~\ref{Thm: B} and Corollary~\ref{thm.with dominated splitting}. If the index of $q$ in Theorem \ref{Thm: B} is no less than $i+k$, then we can conclude Theorem~\ref{Thm: B} directly from Corollary~\ref{thm.with dominated splitting}. If the index of $q$ is smaller than $i+k$, indeed by the no-domination assumption along $F$ and the technics of~\cite{BB}, we can do an arbitrarily small perturbation to get a new hyperbolic periodic orbit of index $i+k$. However, we do not know whether or not the new generated periodic orbits are still contained in the homoclinic class.

The proof of Theorem~\ref{Thm: B} is not by finding a hyperbolic periodic orbit of index $i+k$ in the homoclinic class.  We use a little different strategy from the proof of Theorem~\ref{Thm: A}  to give the proof.

A more general statement than Theorem~\ref{Thm: B} can be expected to be true. We state it as the following question.

\begin{Question}
For generic $f\in\diff^1(M)$, consider the finest dominated splitting $E_1\oplus\cdots\oplus E_k$ over a homoclinic class $H(P,f)$. Assume that there exist two saddles $q_1,q_2$ in the class such that $|\Jac({\rm Df}^{\pi(q_1)}|_{E_i(q_1)})|>1$ and $|\Jac({\rm Df}^{\pi(q_2)}|_{E_j(q_2)})|<1$ where $i\leq j$. Then for any $i\leq l\leq j$, does there exist an ergodic measure whose Lyapunov exponents along the bundle $E_l$ all vanish?
\end{Question}

%%%%%%%%%%%%%%%%%%%%%%%%%%%%%%%%%%%%%%%%%%%%%%%%%%%
\subsection*{Organization of the paper}
%%%%%%%%%%%%%%%%%%%%%%%%%%%%%%%%%%%%%%%%%%%%%%%%%%%%%%

In Section~\ref{Section:preliminary}, we give some definitions and some known results.
Section~\ref{Section:proof of thm A} and Section~\ref{Section:proof of thm B} give the proof of Theorem \ref{Thm: A} and Theorem~\ref{Thm: B} respectively. Section~\ref{Section:counterexample} gives an example  which shows that the assumption of existence of both periodic orbits of index $\dim(E)$ and $\dim(E\oplus F)$ in   Corollary~\ref{thm.with dominated splitting} is important.

%%%%%%%%%%%%%%%%%%%%%%%%%%%%%%%%%%%%%%%%%%%%%%%%%%%%%%%%%%%%%%%%%%%%%%%%%%%%%
%%%%%%%%%%%%%%%%%%%%%%%%%%%%%%%%%%%%%%%%%%%%%%%%%%%%%%%%%%%%%%%%%%%%%%%%%%%%%%%%
\section{Preliminary}\label{Section:preliminary}
%%%%%%%%%%%%%%%%%%%%%%%%%%%%%%%%%%%%%%%%%%%%%%%%%%%%%%%%%%%%%%%%%%%%%%%%%%%%%%%
%%%%%%%%%%%%%%%%%%%%%%%%%%%%%%%%%%%%%%%%%%%%%%%%%%%%%%%%%%%%%%%%%%%%%%%%%%%%%%%%%
In this section, we collect the notations and known results that we need in this paper.

%%%%%%%%%%%%%%%%%%%%%%%%%%%%%%%%%%%%%%%%%%%%%%%%%%%%%%%%%%%%%%%%%%%%%%%
\subsection{Lyapunov exponents}
%%%%%%%%%%%%%%%%%%%%%%%%%%%%%%%%%%%%%%%%%%%%%%%%%%%%%%%%%%%%%%%%%%%%%%%
In this subsection, we state an expression of Lyapunov exponents for an ergodic measure.

Let $f\in\diff^1(M)$ and $\nu$ be an $f$-ergodic   measure. We denote by $$\chi_{1}(\nu,f)\leq\cdots\leq\chi_{d}(\nu,f)$$ all the Lyapunov exponents of $\nu$ counted by multiplicity.
We define a continuous function on $M$ as:
$$L_i^n(x,f)=\frac{1}{n}\log{\norm{\wedge^i {\rm Df}^n(x)}}.$$
Then, for $\nu$-a.e. $x\in M$, we have that $$\chi_{i}(\nu,f)=\lim_{n\rightarrow\infty}(L_{d-i+1}^n(x,f)-L_{d-i}^n(x,f)).$$

%%%%%%%%%%%%%%%%%%%%%%%%%%%%%%%%%%%%%%%%%%%%%%%%%%%%%%%%%%%%%%%%%%%
\subsection{Chain recurrence}
%%%%%%%%%%%%%%%%%%%%%%%%%%%%%%%%%%%%%%%%%%%%%%%%%%%%%%%%%%%%%%%%%%%%%

Let $(X,d)$ be a compact metric space and $f$ be a homeomorphism on $X$. Given two points $x,y\in X$, we define the relation $x\dashv y$, if and only if  for any $\epsilon>0$, there exist finite points $x=z_0,z_1,\cdots,z_k=y$, where $k\geq 1$, such that
$$d(f(z_i),z_{i+1})\leq\epsilon, \textrm{ for any } 0\leq i\leq k-1.$$
We define the relation $x\vdash\!\dashv y$ if and only if $x\dashv y$ and $y\dashv x$.

The \emph{ chain recurrent set} of $f$ is defined as
$$\mathcal{R}(f)=\{x\in X: x\dashv x\}.$$
It's well known that $\vdash\!\dashv$ is an equivalent relation on $\mathcal{R}(f)$. Hence, $\mathcal{R}(f)$ can be decomposed into different equivalent classes, each of which is called a \emph{chain recurrence class}.
\bigskip

Homoclinic classes can also be defined in the following way.

\begin{definition}\label{Def:homoclinic class}
Assume that $f$ is a diffeomorphism in $\diff^1(M)$ and $P,Q$ are two hyperbolic periodic orbits of $f$. We say that $P$ and $Q$ are \emph{homoclinically related}, if $W^u(P)$ has non-empty transverse intersections with $W^s(Q)$, and vice versa, denoted by $W^u(P)\pitchfork W^s(Q))\neq \emptyset$ and $W^s(P)\pitchfork W^u(Q)\neq\emptyset$. We call the closure of the set of periodic orbits homoclinically related to $P$ the \emph{homoclinic class} of $P$ and denote it as $H(P,f)$ or $H(P)$ for simplicity.
\end{definition}

The following lemma is from~\cite{BC}.

\begin{lemma}
For $C^1$-generic diffeomorphisms, the chain recurrence class of a hyperbolic periodic orbit $Q$ coincides with its homoclinic class $H(Q)$.
\end{lemma}

%%%%%%%%%%%%%%%%%%%%%%%%%%%%%%%%%%%%%%%%%%%%%%%%%%%%%%%%%%%%%%%%%%%%%%%%%%%%%
\subsection{A criterion to the ergodicity of convergence}
%%%%%%%%%%%%%%%%%%%%%%%%%%%%%%%%%%%%%%%%%%%%%%%%%%%%%%%%%%%%%%%%%%%%%%%%%%%%%

The period of a periodic orbit $P=\orb(p)$ is denoted by  $\pi(P)$.
We define a relation between two periodic orbits called \emph{good approximation} which is given in~\cite{DG,BDG}.

\begin{definition}\label{Def:good approximation}
 Given a dynamical system $(K,f)$. Let $X,Y$ be two periodic orbits of $f$. We say that $X$ is a \emph{$(\delta, \kappa)$-good approximation} of $Y$, for some $\delta>0$ and $\kappa\in(0,1]$ if there exist a subset $\tilde{X}\subset X$ and a map $\Pi :\tilde{X}\rightarrow Y$ such that:
\begin{itemize}\item $\frac{\#\tilde{X}}{\# X}>\kappa$;
\item $\#(\Pi^{-1}(y))$ is independent of $y$, where $y$ belongs to $Y$;
\item $d(f^{i}(x),f^{i}(\Pi(x)))<\delta$, for any $i=0,\cdots,\pi(Y)-1$ and any $x\in\tilde{X}$.
\end{itemize}
\end{definition}

 %For a sequence of compact subset $K_n$ of a compact metric space, the \emph{topological upper limit} of $K_n$ is defined as:
   %\begin{displaymath}
      %\limsup_{n\rightarrow\infty} K_n=\cap_{n=1}^{\infty}\overline{\cup_{k=n}^{\infty}K_k}.
    %\end{displaymath}

Here, we state a criterion which is first used in \cite{GIKN} and developed in \cite{DG,BDG} showing that with some good approximation assumption, a sequence of periodic measures  converges to an ergodic measure.
\begin{lemma}[\cite{DG,BDG}] \label{limit}
Given a system $(K,f)$. Let $\{X_n\}$ be a sequence of periodic orbits. Assume that $X_{n+1}$ is a $(\delta_n,\kappa_n)$-good approximation of $X_n$ for each $n\in\N$, where $\{\delta_n\}$ and $\{\kappa_n\}$ are two sequences of positive numbers no more than $1$ satisfying:
\begin{displaymath}
 \sum_{n\geq 0}\delta_n<\infty \textrm{ and } \prod_{n\geq 0}\kappa_n\in(0,1].
\end{displaymath}
Then the dirac measure supported on $X_n$ converges to an ergodic measure $\nu$ and the support of $\nu$ is given by
% $\limsup_{n\rightarrow\infty} X_n$.
$$\cap_{n=1}^{\infty}\overline{\cup_{k=n}^{\infty}X_k}.$$
\end{lemma}

%%%%%%%%%%%%%%%%%%%%%%%%%%%%%%%%%%%%%%%%%%%%%%%%%%%%%%%%%%%%%%%%%%%%%%%%%%%%%%%%
\subsection{Perturbation technics}
%%%%%%%%%%%%%%%%%%%%%%%%%%%%%%%%%%%%%%%%%%%%%%%%%%%%%%%%%%%%%%%%%%%%%%%%%%%%%%%%%

  Let $A_1,\cdots,A_l\in GL(d,\mathbb{R})$ and we denote by $B=A_l\comp A_{l-1}\comp\cdots\comp A_1.$  Let $\lambda_{1}(B),\cdots,\lambda_{d}(B)$ be the eigenvalues of $B$,  counted by multiplicity and satisfying $$|\lambda_{1}(B)|\leq\cdots\leq|\lambda_{d}(B)|.$$
  The $i^{th}$ Lyapunov exponent of $B$ is defined as $$\chi_{i}(B)=\frac{1}{l}\log|\lambda_{i}(B)|.$$
We say that $B$ has simple spectrum if all the Lyapunov exponents of $B$ are mutually different.

 We state a  version of Theorem 4.11 in \cite{BB} adapted to our situation. A similar result can be found in~\cite{G}.

\begin{lemma}\label{lemma:bb}
   For any $d\geq 2$, $\epsilon>0$,
and $R>1$, there exist two positive integers $T, l_0$ such that:

  Given $l$ linear maps $A_1,\cdots,A_l\in GL(d,\mathbb{R})$ with $l\geq l_0$ such that $\norm{A_i}, \norm{A_{i}^{-1}}<R$. Assume that $B=A_l\comp A_{l-1}\comp\cdots\comp A_1$ has no $T$-domination of index $j$
  for any $j\in\{i_0+1,\cdots,i_0+k_0-1\}$

    For any $k_0$ numbers $\xi_1,\cdots,\xi_{k_0}$ satisfying:
    \begin{itemize}
    \item[--] $\xi_{k_0}\geq\cdots\geq \xi_{1}$;
    \item[--] $\sum_{i=1}^j\xi_i\geq \sum_{i=1}^j\chi_{i_0+i}(B)$, for any $j=1,\cdots,k_0;$
    \item[--]  $\sum_{i=1}^{k_0}\xi_i= \sum_{i=1}^{k_0}\chi_{i_0+i}(B)$.
    \end{itemize}

   Then there exist $l$ one-parameter families of linear maps $\{(A_{i,t})_{t\in[0,1]}\}_{i=1}^l$ such that:
  \begin{enumerate} \item $A_{i,0}=A_i$ for each $i$;
  \item $\norm{A_{i,t}-A_i}<\epsilon$ and $\norm{A_{i,t}^{-1}-A_{i}^{-1}}<\epsilon$, for each $i$ and any $t\in[0,1]$;
  \item Consider the linear map $B_t=A_{l,t}\comp A_{l-1,t}\comp\cdots\comp  A_{1,t}$, then the Lyapunov exponents of $B_t$ satisfy the following:
      \begin{itemize}
      \item[--] $\chi_{j}(B_t)=\chi_{j}(B)$, for any integer $j\in[1, i_0]\cup[i_0+k_0+1, d]$;
      \item[--] $\sum_{j=1}^{k_0}\chi_{i_0+j}(B_t)=\sum_{j=1}^{k_0}\chi_{i_0+j}(B)$, for any $t\in[0,1]$;
      \item[--] For any $j\in[1, k_0]$, the function $\sum_{i=1}^{j}\chi_{i_0+i}(B_t)$ with respect to variable $t$ is non-decreasing;
      \item[--] $\chi_{i_0+j}(B_1)=\xi_j$, for any $j=1,\cdots,k_0$.
      \end{itemize}
  \end{enumerate}
  \end{lemma}

\begin{remark}\label{r.BB}
In particular, we can take $\xi_1=\cdots=\xi_{k_0}=\frac{1}{k_0}\sum_{i=1}^{k_0}\chi_{i_0+i}(B)$ in Lemma \ref{lemma:bb}.
\end{remark}

The following lemma shows that for periodic orbit of large period, we can do certain small perturbation to make it have simple spectrum.
\begin{lemma}\label{l.get simple spectrum}\cite[Lemma 6.6]{BC}
Given a positive number $K$. For any $\epsilon>0$, there exists an integer $N$ such that
for any $n\geq N$ and any matrices $A_1,\cdots,A_n$ in $GL(2,\mathbb{R})$ satisfying that
$\norm{A_i}<K$ and $\norm{A_i^{-1}}<K$ for any $i=1,\cdots,n$.

Then there exist  matrices $B_1,\cdots, B_n$ in $GL(2,\mathbb{R})$ such that
\begin{itemize}
\item[--] $\norm{A_i-B_i}<\epsilon$ and $\norm{A_i^{-1}-B_i^{-1}}<\epsilon$ for any $i=1,\cdots, n$;
\item[--] the matrix $B_n\comp\cdots\comp B_1$ has simple spectrum.
\end{itemize}
\end{lemma}

\begin{remark}
The original statement of Lemma 6.6 in \cite{BC} is for the matrices in $SL(2,\mathbb{R})$,
but with the assumption that the norm of the matrices and its inverse are uniformly bounded,
the  same conclusion is also true directly from \cite[Lemma 6.6]{BC}.
\end{remark}

We state a generalized Franks lemma by N. Gourmelon~\cite{Go}, which allows us to do a Franks-type perturbation along a hyperbolic periodic orbit which keeps some homoclinic or heteroclinic intersections.
\begin{lemma}[Franks-Gourmelon Lemma]\label{lemma:gourmelon} Given $\epsilon>0$, a diffeomorphism $f\in\diff^1(M)$ and a hyperbolic periodic orbit $Q=\orb(q)$ of period $n$. Consider $n$ one-parameter families of linear maps $\{(A_{i,t})_{t\in[0,1]}\}_{i=0}^{n-1}$ in $GL(d, \mathbb{R})$ satisfying the following properties:
\begin{itemize}\item[--] $A_{i,0}={\rm Df}(f^i(q))$ for any integer $i\in[0, n-1]$;
\item[--] $\norm{A_{i,t}-{\rm Df}(f^i(q))}<\epsilon$ and $\norm{A_{i,t}^{-1}-{\rm Df}^{-1}(f^{i+1}(q))}<\epsilon$, for any $t\in[0,1]$;
\item[--]$A_{n-1,t}\comp\cdots\comp A_{0,t}$ is hyperbolic for any $t\in[0,1]$.
\end{itemize}
Then for any neighborhood $U$ of $Q$, any number $\eta>0$ and any pair of compact sets $K^s\subset W_{\eta}^s(Q,f)$ and $K^{u}\subset W_{\eta}^{u}(Q,f) $ which do not intersect $U$, there is a diffeomorphism $g\in\diff^1(M)$ which is $\varepsilon$-$C^1$-close to $f$, such that

\begin{itemize}
\item $g$ coincides with $f$ on $Q\cup M\backslash U$,
\item ${\rm Dg}(g^i(q))=A_{i,1}$, for any $i\in\{0,1,\cdots,n-1\}$,
\item $K^s\subset W^s(Q,g)$ and $K^u\subset W^u(Q,g)$.
\end{itemize}
\end{lemma}

\begin{definition}
Consider a diffeomorphism $f\in\diff^1(M)$. An invariant compact set $\Lambda$ is said to admit a \emph{partially hyperbolic splitting}, if there is a splitting $T_{\Lambda}M=E^s\oplus E^c\oplus E^u$ such that, the splittings $(E^s\oplus E^c)\oplus E^u$ and $E^s\oplus (E^c\oplus E^u)$ are dominated splittings, and the  bundle $E^s$ (resp. $E^u$) is uniformly contracting (resp. expanding). Moreover, at least one of the two extreme bundles $E^s$ and $E^u$ is non-degenerate.
\end{definition}

Consider two hyperbolic periodic points $p$ and $q$ of indices $i$ and $i+k$ respectively.
We say that $p$ and $q$ form a \emph{heterodimensional cycle}
if $W^u(P)$ has   transverse intersections with  $W^s(Q)$ along the orbit of some point $y$, and $W^s(P)$ has   quasi-transverse intersections with  $W^u(Q)$ along the orbit of some point $x$, i.e. $T_x W^s(P)+T_x W^u(Q)$ is a direct sum.
 We say $p$ and $q$ form a \emph{partially hyperbolic heterodimensional cycle},
 if the $f$-invariant compact set $\mathcal{C}=\overline{\orb(x)}\cup \overline{\orb(y)}$ admits a partially hyperbolic
 splitting of the form $T_{\mathcal{C}}M=E^s\oplus E^c\oplus E^{u}$, where $\dim(E^s)=i$ and $\dim(E^c)=k$.
 %, and the two extreme  bundles $E^s$ and $E^u$ are uniformly contracted and expanded by ${\rm Df}$ respectively.
Moreover, for any $x\in\mathcal{C}$, we denote by $W^{ss}(x)$ (resp. $W^{uu}(x)$)
 the \emph{strong stable manifold} (resp. \emph{strong unstable manifold}) of $x$
 which is tangent to the bundle $E^s$ (resp. $E^u$) at $x$.

We have the following theorem from~\cite{BDPR} to obtain \emph{transition between two periodic orbits of different indices}.

\begin{theorem}$\cite[\textrm{Theorem  3.1  and   Lemma 3.5}]{BDPR}$\label{Thm:transition}
Consider a diffeomorphism $f\in\diff^1(M)$. Let $p$ and $q$ be two hyperbolic periodic points of indices $i$ and $i+k$ respectively and denote by $P$ and $Q$ their orbits respectively. Assume that there exist dominated splitting  $T_{P}M=E_1(P)\oplus E_2(P)\oplus E_3(P)$ and $T_{Q}M=E_1(Q)\oplus E_2(Q)\oplus E_3(Q)$ satisfying that $\dim(E_1(P))=\dim(E_1(Q))=i$ and $\dim(E_2(P))=\dim(E_2(Q))=k$. Assume, in addition, that $P$ and $Q$ form a heterodimensional cycle.
Denote by $M_P$ and $M_Q$ the two linear maps: $${\rm Df}^{\pi(P)}(p):T_pM\rightarrow T_pM \text{ and } {\rm Df}^{\pi(Q)}(q):T_qM\rightarrow T_qM.$$

Then for any $C^1$-neighborhood $\mathcal{U}$ of $f$, for any two neighborhoods $U_P$ and $U_Q$ of $P$ and $Q$ respectively, there are two matrices $T_0$ and $T_1$, and two integers $t_0$ and $t_1$, such that  for any two positive integers $m$ and $n$, there is a diffeomorphism $g\in\mathcal{U}$  with a periodic point $p_1$, satisfying the following properties:
\begin{itemize}
\item[--] $g$ and ${\rm Dg}$ coincide with $f$ and ${\rm Df}$ on $P\cup Q$  respectively;
\item[--] For $i=1,2,3$, we have that $T_0(E_i(p))=E_i(q)$ and $T_1(E_i(q))=E_i(p)$;
\item[--] The period of $p_1$ equals $t_0+t_1+n\pi(P)+m\pi(Q)$;
\item[--] the matrix ${\rm Dg}^{\pi(p_1)}(p_1):T_{p_1}M\rightarrow T_{p_1}M$ is conjugate to $$T_1\circ M_Q^m \circ T_0\circ M_P^n;$$
\item[--] we denote by $P_1$ the orbit of $p_1$ under $g$, then we have:
$$\#(P_1\cap U_P)\geq n\pi(P) \text{ and } \#(P_1\cap U_Q)\geq m\pi(Q).$$
\end{itemize}

\end{theorem}
\begin{remark}\label{r.transition}
By Lemma 4.13 in \cite{BDP}, if the periodic orbits  $P$ and $Q$ admit another dominated splitting of the same index,
 the two matrices $T_0$ and $T_1$ can be chosen to preserve the two dominated splitting at the same time.
\end{remark}

%%%%%%%%%%%%%%%%%%%%%%%%%%%%%%%%%%%%%%%%%%%%%%%%%%%%%%%%%%%%%%%%%%%%%%%%%%%%%
\subsection{Generic diffeomorphisms}
%%%%%%%%%%%%%%%%%%%%%%%%%%%%%%%%%%%%%%%%%%%%%%%%%%%%%%%%%%%%%%%%%%%%%%%%%%%%%%
Let $f\in\diff^1(M)$, $P$ and $Q$ be two hyperbolic periodic orbits of $f$. We say that $P$ and $Q$ are \emph{robustly in the same chain recurrence class}, if there exists a $C^1$ small neighborhood $\mathcal{U}$ of $f$ such that for any $g\in\mathcal{U}$, the continuation $P_g$ of $P$ and the continuation $Q_g$ of $Q$ are in the same chain recurrence class.  A periodic orbit $P$ is said to have \emph{simple spectrum}, if the $d$ Lyapunov exponents of $P$ are mutually different.
Denote by $Per(f)$ the set of periodic points of $f$.

The following theorem summarizes some generic properties for $\diff^1(M)$, see for example~\cite{ABCDW,BC,BDPR,BDV,CCGWY,DG}.
 \begin{theorem}\label{Thm:generic properties}
  There exists a residual subset $\mathcal{R}$ of $\diff^1(M)$ such that for any $f\in\mathcal{R}$, we have the followings:
 \begin{enumerate}
 \item\label{generic:K-S}  $f$ is \emph{Kupka-Smale}.

 \item\label{generic:homoclinic class coinsides} Any chain recurrence class containing a hyperbolic periodic orbit $P$ coincides with the homoclinic class $H(P,f)$. Hence two homoclinic classes either coincide or are disjoint.

 \item\label{generic:continuaition of homoclinic class} Given a hyperbolic periodic orbit $P$, there exists a neighborhood $\mathcal{U}$ of $f$ such that the map $g\mapsto H(P_{g},g)$ is well defined and $f$ is a continuous point of this map.
 %\item\label{generic:index complete} For any homoclinic class $H(p,f)$, there exist a $C^1$-neighborhood $\mathcal{V}$ of $f$
 %and  an interval $[a,b]$ of positive integers, such that, for any $g\in\mathcal{V}$, the set of indices of hyperbolic periodic
 % points contained in $H(p_g,g)$ is the interval $[a,b]$.

 \item\label{generic:robust in a chain class}  Given a homoclinic class $H(P,f)$, for any hyperbolic periodic orbit $Q$ contained in $H(P,f)$,
  we have that $P$ and $Q$ are robustly in the same chain recurrence class.
     %there exists a small neighborhood $\mathcal{U}_f$ of $f$ such that for any $g\in\mathcal{U}_f$, the periodic orbits
     %$P_g$ and $Q_g$ are in the same chain recurrence class.

 %\item\label{generic:non-domination} For any homoclinic class $H(P,f)$ admitting a splitting of the form
 %$T_{H(P,f)}M=E^s\oplus E^c\oplus E^u$, if $E^c$ has no finer dominated splitting, then there exists a $C^1$-neighborhood
 %$\mathcal{U}$ of $f$ such that, for any $g\in\mathcal{U}\cap\mathcal{R}$, we have that $T_{H(P_g,g)}M=E^s\oplus E^c\oplus E^u$
 %and $E^c$ has no finer domination, where $P_g$ is the continuation of $P$.

 \item\label{generic:periodic orbits with simple spectrum} Consider a non-trivial homoclinic class $H(P,f)$,
 the set
 $$\{q\in Per(f): \text{ $\orb(q)$ has simple spectrum and is homoclinically related to $P$}\}$$
 is dense in $H(P,f)$.

 \item\label{generic:get cycle}Consider a hyperbolic periodic orbit $P$ of index $i$, whose homoclinic class contains  a hyperbolic periodic orbit  $Q$ of index $i+k$ for some integer $k>0$. If there exist dominated splitting  $T_{P}M=E_1(P)\oplus E_2(P)\oplus E_3(P)$ and $T_{Q}M=E_1(Q)\oplus E_2(Q)\oplus E_3(Q)$ satisfying that $\dim(E_1(P))=\dim(E_1(Q))=i$ and $\dim(E_2(P))=\dim(E_2(Q))=k$.  Then arbitrarily $C^1$-close to $f$, there is a diffeomorphism $g$, satisfying that:
     \begin{itemize}
     \item[--] $g$ and ${\rm Dg}$ coincide with $f$ and ${\rm Df}$ on $P\cup Q$ respectively,
     \item[--] under the diffeomorphism $g$, the periodic orbits $P$ and $Q$ form a partially hyperbolic heterodimensional cycle $K$, which is contained in an arbitrarily small neighborhood of $H(P,f)$.
     \end{itemize}

 \item\label{generic:shadowing periodic orbits} Consider a hyperbolic periodic orbit $P$ with simple spectrum whose homoclinic class $H(P,f)$ is non-trivial. Then for any $\epsilon>0,\delta>0$ and $\kappa\in(0,1)$, there is a hyperbolic periodic orbit $Q$ homoclinically related to $P$, such that the following properties are satisfied:
    \begin{itemize}
    \item[--] $Q$ has simple spectrum and is $\varepsilon$-dense in $H(P,f)$;
         \item[--] $Q$ is a $(\delta,\kappa)$-good approximation of $P$;
     \item[--] $|\chi_{i}(Q,f)-\chi_{i}(P,f)|<\epsilon$, for any $i\in\{1,2,\cdots,d\}$.
    \end{itemize}
 \end{enumerate}
 \end{theorem}

%%%%%%%%%%%%%%%%%%%%%%%%%%%%%%%%%%%%%%%%%%%%%%%%%%%%%%%%%%%%%%%%%%%%%%%%%%%%
%%%%%%%%%%%%%%%%%%%%%%%%%%%%%%%%%%%%%%%%%%%%%%%%%%%%%%%%%%%%%%%%%%%%%%%%%%%%%%
\section{Ergodic measure with multi-zero   Lyapunov exponents for the case  controlled by norm: Proof of Theorem~\ref{Thm: A}}\label{Section:proof of thm A}
%%%%%%%%%%%%%%%%%%%%%%%%%%%%%%%%%%%%%%%%%%%%%%%%%%%%%%%%%%%%%%%%%%%%%%%%%%%%%%%%
%%%%%%%%%%%%%%%%%%%%%%%%%%%%%%%%%%%%%%%%%%%%%%%%%%%%%%%%%%%%%%%%%%%%%%%%%%%%%%%%%

%%%%%%%%%%%%%%%%%%%%%%%%%%%%%%%%%%%%%%%%%%%%%%%%%%%%%%%%%%%%%%%%%%%%%%%%%%%%%%
  \subsection{ Proof of Theorem~\ref{Thm: A}}
%%%%%%%%%%%%%%%%%%%%%%%%%%%%%%%%%%%%%%%%%%%%%%%%%%%%%%%%%%%%%%%%%%%%%%%%%%%%%%%%
  The following proposition is the main step for proving Theorem~\ref{Thm: A}.

  \begin{proposition}\label{proposition}
For generic $f\in\diff^1(M)$. Consider a non-trivial homoclinic class $H(P,f)$ of a hyperbolic periodic orbit $P$ of index $i$. Assume that
\begin{itemize}
\item[--] there is a hyperbolic  periodic orbit $Q$ of index $i+k$ contained in $H(P,f)$, where $k\geq 1$;
\item[--] there is no dominated splitting of index $i+j$ over $H(P,f)$, for any $j=1,2,\cdots,k-1$.
\end{itemize}
Then there is a constant $\chi>0$  such that for any $\gamma>0$ and   any hyperbolic periodic orbit  $P_0$ with simple spectrum which is homoclinically related to $P$,  there is a hyperbolic periodic orbit $P_1$ homoclinically related to $P$, satisfying the following properties:
\begin{enumerate}
\item $P_1$ is $\gamma$-dense in $H(P,f)$ and has simple spectrum;
\item $\chi_{i+k}(P_1,f)<\frac{3}{4}\cdot\chi_{i+k}(P_{0},f)$;
\item $P_1$ is a $\big(\gamma, 1-\frac{\chi_{i+k}(P_0,f)}{\chi+\chi_{i+k}(P_0,f)}\big)$-good approximation of $P_0$.
\end{enumerate}
\end{proposition}
\vspace{2mm}
Using Proposition~\ref{proposition}, we give the proof of Theorem~\ref{Thm: A}.

 \proof[Proof of Theorem~\ref{Thm: A}]
 By item~\ref{generic:homoclinic class coinsides} of Theorem \ref{Thm:generic properties}, we can assume that $P$ is of index $i$. We take the positive constant $\chi$ from Proposition~\ref{proposition}.
 We will inductively construct a sequence of periodic orbits $\{P_n\}$ with simple spectrum, a sequence of positive numbers $\{\gamma_n\}$ and a sequence of integers $\{N_n\}$ satisfying the following properties:
\begin{enumerate}
 \item $\chi_{i+k}(P_{n+1},f)<\frac{3}{4}\cdot\chi_{i+k}(P_n,f)$;
 \item $P_n$ is homoclinically related to $P$ and is $\frac{1}{2^n}$-dense inside $H(P,f)$;
 \item the constants $\gamma_{n} $ and $N_{n} $  satisfy that:
   \begin{itemize}
   \item[--]  $\gamma_{n}<\frac{1}{2}\gamma_{n-1}$  and $ N_{n}>N_{n-1};$
   \item[--] for any point $x\in B_{2\gamma_{n}}(P_n)\cap H(P,f)$, we have that
      $$0<L_{d-i}^{N_n}(x,f)-L_{d-i-1}^{N_n}(x,f),$$
     $$0<L_{d-i}^{N_n}(x,f)-L_{d-i-k}^{N_n}(x,f)<2k\cdot\chi_{i+k}(P_n,f);$$
     \end{itemize}
 \item $P_{n+1}$ is $\big(\gamma_{n}, 1-\frac{\chi_{i+k}(P_n,f)}{\chi+\chi_{i+k}(P_n,f)}\big)$-good approximation of $P_n$.

     %\begin{enumerate} \item  $\frac{1}{N_{n}}\log\norm{{\rm Df}^{N_{n}}|_{E^c(x)}}<2\chi_{i+k}(P_n)$;
     %\item $\frac{1}{N_{n}}\log m({\rm Df}^{N_{n}}|_{E^c(x)})>-2\chi_{i+k}(P_n)$.
     %\end{enumerate}
\end{enumerate}

\paragraph{Choice of $P_0$, $N_0$ and $\gamma_0$}

First we construct for $n=0$. By the item \ref{generic:periodic orbits with simple spectrum} and item~\ref{generic:shadowing periodic orbits} of Theorem \ref{Thm:generic properties}, we can  choose a hyperbolic periodic orbit $P_0$, with simple spectrum, which is homoclinically related to $P$ and is $\frac{1}{2^0}$-dense inside $H(P,f)$. Hence the item $2$ is satisfied.

By the definition of the function $L_j^n(x,f)$,  there exists an integer $N_{0}$  such that
for any $y\in P_{0}$, we have that
   $$0<L_{d-i}^{N_0}(y,f)-L_{d-i-1}^{N_0}(y,f),$$
  $$0<L_{d-i}^{N_0}(y,f)-L_{d-i-k}^{N_0}(y,f)<\frac{3k}{2}\cdot\chi_{i+k}(P_0,f).$$
By the uniform continuity of the functions $L_{d-i}^{N_0}(x,f)$ and $L_{d-i-k}^{N_0}(x,f)$, there exists  a number $\gamma_{0}>0$ such that
for any point $x\in B_{2\gamma_{0}}(P_0)\cap H(P,f)$, we have that
   $$0<L_{d-i}^{N_0}(x,f)-L_{d-i-1}^{N_0}(x,f),$$
  $$0<L_{d-i}^{N_0}(x,f)-L_{d-i-k}^{N_0}(x,f)<2k\cdot\chi_{i+k}(P_0,f).$$
Hence the item $3$ is satisfied. Notice that we do not have to check the items $1,4$ for $n=0$.
\medskip

\paragraph{Construct $P_n$, $N_n$ and $\gamma_n$ inductively}
Assume that $P_j$, $N_j$ and $\gamma_j$ are already defined for any $j\leq n$.
We apply $P_n$, $\gamma_n$, and $\frac{1}{2^{n+1}}$ to Proposition \ref{proposition}, then we get a periodic orbit $P_{n+1}$ with simple spectrum, satisfying that:
\begin{itemize}
 \item $\chi_{i+k}(P_{n+1},f)<\frac{3}{4}\cdot\chi_{i+k}(P_n,f)$;
 \item $P_{n+1}$ is homoclinically related to $P$ and is $\frac{1}{2^{n+1}}$-dense in $H(P,f)$;
 \item $P_{n+1}$ is $(\gamma_{n},1-\frac{\chi_{i+k}(P_n,f)}{\chi+\chi_{i+k}(P_n,f)})$-good approximation of $P_n$.
\end{itemize}
Then the items $1,2,4$ are satisfied.

By the definition of the function $L_j^n(x,f)$, there is an integer $N_{n+1}>N_n$  satisfying that:
  for any $y\in P_{n+1}$, we have that
  $$0<L_{d-i}^{N_{n+1}}(y,f)-L_{d-i-1}^{N_{n+1}}(y,f),$$
  $$0<L_{d-i}^{N_{n+1}}(y,f)-L_{d-i-k}^{N_{n+1}}(y,f)<\frac{3k}{2}\cdot\chi_{i+k}(P_{n+1},f).$$
 By the uniform continuity of the functions   $L_{d-i}^{N_{n+1}}(x,f)$ and $L_{d-i-k}^{N_{n+1}}(x,f)$, there exists  a number $\gamma_{n+1}\in(0,\frac{1}{2}\gamma_n)$ such that
  for any point $x\in B_{2\gamma_{n+1}}(P_{n+1})\cap H(P,f)$, we have
 $$0<L_{d-i}^{N_{n+1}}(x,f)-L_{d-i-1}^{N_{n+1}}(x,f),$$
 $$0<L_{d-i}^{N_{n+1}}(x,f)-L_{d-i-k}^{N_{n+1}}(x,f)<2k\cdot\chi_{i+k}(P_{n+1},f).$$

\paragraph{End of proof of Theorem~\ref{Thm: A}}
 By Lemma~\ref{limit}, the sequence of ergodic measures $\delta_{P_n}$ converges to an ergodic measure $\nu$ whose support is $H(P,f)$. 
 We will show that $\nu$ has $k$ vanishing Lyapunov exponents. %all the center Lyapunov exponents of $\nu$ vanish.

 \begin{claim}\label{c.zero Lyapunov exponent}
 The $(i+j)^{th}$ Lyapunov exponent of $\nu$ equals zero, for any $j=1,2,\cdots,k$.
 \end{claim}

\proof
By Definition~\ref{Def:good approximation}, there exist   a subset $\tilde{P_n}$ of $P_n$ and a map $\Pi_n :\tilde{P_n}\mapsto P_{n-1}$ for each $n\geq 2$.
Consider the set $K_n=\Pi_{n}^{-1}\comp\Pi_{n-1}^{-1}\comp\cdots\comp \Pi_{1}^{-1}(P_0)$,
then we have that $$\delta_{P_n}(K_n)\geq\prod_{l=0}^{n-1}\big(1-\frac{\chi_{i+k}(P_l,f)}{\chi+\chi_{i+k}(P_l,f)}\big).$$
We denote by  $$K=\cap_{n=1}^{\infty}\overline{\cup_{l=n}^{\infty}K_l} ,$$
then we have that $\nu(K)\geq \lim_{n\rightarrow\infty}\delta_{P_n}(K_n)>0$.

On the other hand, for any point $x\in B_{2\gamma_{n}}(P_n)\cap H(P,f)$, we have that
     $$0<L_{d-i}^{N_n}(x,f)-L_{d-i-1}^{N_n}(x,f),$$
     $$0<L_{d-i}^{N_n}(x,f)-L_{d-i-k}^{N_n}(x,f)<2k\cdot\chi_{i+k}(P_n,f).$$

Since $P_{n+1}$ is a $(\gamma_{n},1-\frac{\chi_{i+k}(P_n)}{\chi+\chi_{i+k}(P_n)})$-good approximation of $P_n$, we have that $K$ is contained in the $\sum_{i=n}^{\infty}\gamma_i$ neighborhood of $P_n$, therefore is contained in $2\gamma_{n}$ neighborhood of $P_n$.
As a consequence,  for any $y\in K$, we have the following
  \begin{equation}\label{e.positive center exponent}
  0< L_{d-i}^{N_n}(y,f)-L_{d-i-1}^{N_n}(y,f),
   \end{equation}
   \begin{equation}\label{e.estimate of exponents}
   0<L_{d-i}^{N_n}(y,f)-L_{d-i-k}^{N_n}(y,f)<2k\cdot\frac{3^n}{4^{n}}\cdot\chi_{i+k}(P,f).
   \end{equation}

Since $\nu$ is ergodic,  for $\nu$-a.e. point $y$,  we have that
   \begin{equation}\label{e.sum of center exponents}
    \sum_{j=1}^{k}\chi_{i+j}(\nu,f)= \lim_{n\rightarrow+\infty}\big( L_{d-i}^{n}(y,f)-L_{d-i-k}^{n}(y,f)\big),
   \end{equation}
   \begin{equation}\label{e.i th exponent}
    \chi_{i+1}(\nu,f)=\lim_{n\rightarrow+\infty}\big(L_{d-i}^{n}(y,f)-L_{d-i-1}^{n}(y,f)\big).
   \end{equation}
By the fact that $\nu(K)>0$ and the formulas~(\ref{e.estimate of exponents}) and~(\ref{e.sum of center exponents}), we can see that $$\sum_{j=1}^{k}\chi_{i+j}(\nu,f)=0.$$
By the formulas~(\ref{e.positive center exponent}) and~(\ref{e.i th exponent}), we get that $\chi_{i+1}(\nu,f)\geq 0$.
Then by the fact that $\chi_{i+1}(\nu,f)\leq \chi_{i+2}(\nu,f)\leq\cdots\leq \chi_{i+k}(\nu,f)$,  we have that
$$\chi_{i+j}(\nu,f)= 0, \textrm{ for any $j=1,2,\cdots,k$.}$$

This ends the proof of Theorem~\ref{Thm: A}.
 \endproof
Now it remains to prove Proposition~\ref{proposition}.

%%%%%%%%%%%%%%%%%%%%%%%%%%%%%%%%%%%%%%%%%%%%%%%%%%%%%%%%%%%%%%%
\subsection{Good approximation with weaker center Lyapunov exponents: Proof of Proposition~\ref{proposition}}
%%%%%%%%%%%%%%%%%%%%%%%%%%%%%%%%%%%%%%%%%%%%%%%%%%%%%%%%%%%%%%%

The proof of Proposition~\ref{proposition} is based on the following perturbation Lemma:
\begin{lemma}\label{main lemma}
Consider a diffeomorphism $f\in\diff^1(M)$.
Let  $P$ and $Q$ be two hyperbolic periodic orbits of indices $i$ and $i+k$ respectively. Assume that  $Q$ and $P$ form a partially hyperbolic heterodimensional cycle $K$ with the splitting $T_{K}M=E^s\oplus E^c\oplus E^{u}$. Assume, in addition, that
$$\chi_{i+j}(P,f)=\log{\mu}>0 \textrm{  and } \chi_{i+j}(Q,f)=\log{\lambda}<0,  \textrm{ for any $j=1,2,\cdots,k$.} \textrm{ $(\bigstar)$ }$$

Then for any $\gamma>0$ and any $C^1$ neighborhood $\mathcal{U}$ of $f$, there exist
a diffeomorphism $g\in\mathcal{U}$ and a hyperbolic periodic orbit $P_1$ of $g$ such that
 \begin{enumerate}
 \item $P_1$ has simple spectrum;
 \item $g$ and ${\rm Dg}$ coincide with $f$ and ${\rm Df}$  on the set $P\cup Q$ respectively ;
 \item $\frac{1}{4}\cdot\chi_{i+k}(P,g)<\chi_{i+1}(P_1,g)<\chi_{i+k}(P_1,g)<\frac{1}{2}\cdot\chi_{i+k}(P,g)$;
 \item $P_1$ is a $\big(\gamma, 1+\frac{\log\mu}{2\log\lambda-\log\mu}\big)$-good approximation of $P$;
 \item $W^{ss}(P_1)$ has transverse intersections with  $W^{u}(P)$   and $W^{uu}(P_1)$ has transverse intersections with  $W^s (Q)$.
 \end{enumerate}
 \end{lemma}
 \begin{remark}\label{r.robustly in the chain class} If $P$ and $Q$ are robustly in the same chain recurrence class, the last item of Lemma \ref{main lemma} implies that $P_1$ is robustly in the same chain recurrence class with $P$ and $Q$.
 \end{remark}

The idea of the proof of Lemma~\ref{main lemma} is that we  mix two hyperbolic periodic orbits of different indices to get a new periodic orbit with weaker center Lyapunov exponents.

\proof[Proof of Lemma~\ref{main lemma}]   We fix a small number $\gamma>0$ and a neighborhood $\mathcal{U}$ of $f$. There exists $\epsilon>0$ such that the $\epsilon$ neighborhood of $f$ is contained in $\mathcal{U}$. There is a small number $0<\theta<1$, such that for any $h\in\mathcal{U}$ and any two points $z_1,z_2$ satisfying  $d(z_1,z_2)<\theta\cdot\gamma$, we have that $$d(h^i(z_1),h^i(z_2))<\frac{\gamma}{2}, \textrm{  for any $i\in [-\pi(P),\pi(P)]$}.$$
We take two neighborhoods $U_P$ and $U_Q$ of $P$ and $Q$ respectively, such that $U_P$ is contained in the $\theta\cdot\gamma$-neighborhood of $P$ and is disjoint from $U_Q$.

\paragraph{Construction of the periodic orbit $P_1$}
Let  $P=\orb(p)$ and $Q=\orb(q)$. We denote by $M_P$ and $M_Q$ the two linear maps: $${\rm Df}^{\pi(P)}(p):T_pM\rightarrow T_pM \text{ and } {\rm Df}^{\pi(Q)}(q):T_qM\rightarrow T_qM.$$
Since $P$ and $Q$ form a partially hyperbolic heterodimensional cycle $K$, by Theorem~\ref{Thm:transition}, there are two matrices $T_0, T_1$ and two integers $t_0, t_1$ such that for any two integers $m$ and $n$, there is a diffeomorphism $g$,  which is $\frac{\epsilon}{4}$-$C^1$-close to $f$ and has a periodic orbits $P_1=\orb(p_1,g)$, satisfying the following properties:
\begin{itemize}
\item[--] $g$ and ${\rm Dg}$ coincide with $f$ and ${\rm Df}$ respectively on $P\cup Q$,
\item[--] the matrix ${\rm Dg}^{\pi(p_1)}(p_1):T_{p_1}M\rightarrow T_{p_1}M$ is conjugate to $$T_1\circ M_Q^m \circ T_0\circ M_P^n.$$
\item[--] $\pi(p_1)=t_0+t_1+n\pi(P)+m\pi(Q)$.
\item[--] $\#(P_1\cap U_P)\geq n\pi(P), \text{ and } \#(P_1\cap U_Q)\geq m\pi(Q).$
\end{itemize}

Moreover, by the continuity of partial hyperbolicity and the local stable and unstable manifolds of hyperbolic periodic orbit,
by taking $U_P$ and $U_Q$ small enough at first, we have that $W^{ss}(P_1,g)$ intersects $W^{u}_{loc}(P,g)$ transversely and $W^{uu}(P_1,g)$ intersects $W^s_{loc}(Q,g)$ transversely.
\bigskip

By the second item of Theorem \ref{Thm:transition}, we can take    proper coordinates  at $T_PM$ and $T_QM$, under which we have:

\begin{displaymath}
\mathbf{M_P}=
\left(\begin{array}{ccc}
A_{s}&0&0\\
0&A_{c}&0\\
0&0&A_{u}
\end{array}\right),\text{ }
\mathbf{M_Q}=
\left(\begin{array}{ccc}
B_{s}&0&0\\
0&B_{c}&0\\
0&0&B_{u}
\end{array}\right)
\end{displaymath}
\begin{displaymath}
\mathbf{T_1}=
\left(\begin{array}{ccc}
C_{s}&0&0\\
0&C_{c}&0\\
0&0&C_{u}
\end{array}\right),\text{ }
\mathbf{T_0}=
\left(\begin{array}{ccc}
D_{s}&0&0\\
0&D_{c}&0\\
0&0&D_{u}
\end{array}\right)
\end{displaymath}

\paragraph{Choice of the integers $m$ and $n$}
 We will adjust $m,n$ to get the periodic orbit that satisfies the properties stated in Lemma~\ref{main lemma}. We take $\eta>0$ which will be decided later.
 \begin{claim}\label{c.range of center Lyapunov exponent}There exists an integer $N_{\eta}$ such that for any $m\geq N_{\eta}$ and $n\geq N_{\eta}$, we have that
  all the center Lyapunov exponents of $P_1$ belong to  the interval:
$$\Big[\frac{m\cdot\pi(P)\cdot\log\mu+n\cdot\pi(Q)\cdot\log\lambda}{m\cdot\pi(P)+n\cdot\pi(Q)}-2\eta,\,
\frac{m\cdot\pi(P)\cdot\log\mu+n\cdot\pi(Q)\cdot\log\lambda}{m\cdot\pi(P)+n\cdot\pi(Q)}+2\eta\Big].$$
 \end{claim}
 \proof
  By the Equation $(\bigstar)$  in the   assumption of Lemma \ref{main lemma},
   there exists an integer $N_1(\eta)$ such that for any $m, n\geq N_1(\eta)$, we have that
 \begin{align*}
 &\log\mu-\eta<\frac{1}{m\cdot\pi(P)}\log{\um(A_{c}^{m})}\leq\frac{1}{m\cdot\pi(P)}\log{\norm{A_{c}^{m}}}< \log\mu+\eta;
 \\
 &\log\lambda-\eta<\frac{1}{n\cdot\pi(Q)}\log{\um(B_{c}^{n})}\leq\frac{1}{n\cdot\pi(Q)}\log{\norm{B_{c}^{n}}}< \log\lambda+\eta.
 \end{align*}
  As a consequence,  for any unit vector $v\in E^{c}(P_{1})$ and $k\in\N$, we have that
\begin{align*}
  &\norm{{\rm Dg}^{k\cdot\pi(P_1)}v}\leq (\norm{C_{c}}\cdot\norm{D_{c}})^{k}\cdot\exp\big(k\cdot m\cdot \pi(P)\cdot(\log\mu+\eta)+k\cdot n\cdot \pi(Q)\cdot(\log\lambda+\eta)\big),
 \\
  &\norm{{\rm Dg}^{k\cdot \pi(P_1)}v}\geq
  (\um(C_{c})\cdot \um(D_{c}))^{k}\cdot\exp\big(k\cdot m\cdot \pi(P)\cdot(\log\mu-\eta)+k\cdot n\cdot\pi(Q)\cdot(\log\lambda-\eta)\big).
\end{align*}
  Hence, \begin{align*}
  &\frac{1}{k\cdot \pi(P_1)}\log\norm{{\rm Dg}^{k\cdot\pi(P_1)}v}\leq \frac{\log{(\norm{C_{c}}\cdot\norm{D_{c}})}}{\pi(P_1)}+\frac{m\cdot\pi(P)\cdot(\log\mu+\eta)+n\cdot\pi(Q)\cdot(\log\lambda+\eta)}{\pi(P_1)}.
  \\
  &\frac{1}{k\cdot\pi(P_1)}\log\norm{{\rm Dg}^{k\pi(P_1)}v}\geq \frac{\log{(m(C_{c})\cdot m(D_{c}))}}{\pi(P_1)}+\frac{m\cdot\pi(P)\cdot(\log\mu-\eta)+n\cdot\pi(Q)\cdot(\log\lambda-\eta)}{\pi(P_1)}.
\end{align*}

By the fact that $\pi(P_1)=m\pi(P)+n\pi(Q)+t_0+t_1$  and the matrices  $C_c, D_c$ are independent of $m$ and $ n$, there exists  an integer $N_2(\eta)$ such that
for any $m,n\geq N_2(\eta)$, we have that
\begin{itemize}
\item 
 $$-\frac{\eta}{2}<\frac{\log({\um(C_{c})\cdot \um(D_{c})})}{\pi(P_1)}\leq \frac{\log{(\norm{C_{c}}\cdot\norm{D_{c}}})}{\pi(P_1)}<\frac{\eta}{2};$$
\item
$$\Big|\frac{m\cdot\pi(P)\cdot\log\mu+n\cdot\pi(Q)\cdot\log\lambda}{\pi(P_1)}
-\frac{m\cdot\pi(P)\cdot\log\mu+n\cdot\pi(Q)\cdot\log\lambda}{m\pi(P)+n\pi(Q)}\Big|<\frac{\eta}{2}.$$
\end{itemize}
We take $N_{\eta}=\max\{N_1(\eta), N_2(\eta)\}$. When $m,n\geq N_{\eta}$,
we have that  all the center Lyapunov exponents of $P_1$  would belong to  the interval:
$$\Big[\frac{m\cdot\pi(P)\cdot\log\mu+n\cdot\pi(Q)\cdot\log\lambda}{m\cdot\pi(P)+n\cdot\pi(Q)}-2\eta,\,
\frac{m\cdot\pi(P)\cdot\log\mu+n\cdot\pi(Q)\cdot\log\lambda}{m\cdot\pi(P)+n\cdot\pi(Q)}+2\eta\Big].$$
This ends the proof of Claim~\ref{c.range of center Lyapunov exponent}.
\endproof

To guarantee the item $3$,
we only need that
\begin{equation}\label{equa:upper bound}
\frac{m\pi(P)\cdot\log\mu+n\pi(Q)\cdot\log\lambda}{m\pi(P)+n\pi(Q)}+2\eta<\frac{1}{2}\log\mu
\end{equation}
and
\begin{equation}\label{equa:lower bound}
\frac{m\pi(P)\cdot\log\mu+n\pi(Q)\cdot\log\lambda}{m\pi(P)+n\pi(Q)}-2\eta>\frac{1}{4}\log\mu.
\end{equation}

By the choice of the numbers $\theta$ and $\gamma$, to guarantee the item $4$, we only need that
\begin{equation}\label{equa:good approximation}
\frac{m\pi(P)}{m\pi(P)+n\pi(Q)}>1+\frac{\chi_{i+k}(P,g)}{2\chi_{i+k}(Q,g)-\chi_{i+k}(P,g)}=1+\frac{\log\mu}{2\log\lambda-\log\mu}.
\end{equation}

By calculation, to satisfy the inequalities~(\ref{equa:upper bound}),~(\ref{equa:lower bound}) and~(\ref{equa:good approximation}), we only have to show that there exist  $m,n$ large enough such that the following is satisfied:
\begin{equation}\label{equa:sufficient condition}
\max\big\{\frac{-2\log\lambda}{\log\mu},\frac{\log\mu-4\log\lambda+8\eta}{3\log\mu-8\eta}\big\}<\frac{m\pi(P)}{n\pi(Q)}
<\frac{\log\mu-2\log\lambda-4\eta}{\log\mu+4\eta}.
\end{equation}

When $\eta$ is chosen small, we have the following  inequality
\begin{equation}\label{equa:more sufficient condition}
\max\big\{\frac{-2\log\lambda}{\log\mu},\frac{\log\mu-4\log\lambda+8\eta}{3\log\mu-8\eta}\big\}
<\frac{\log\mu-2\log\lambda-4\eta}{\log\mu+4\eta}.
\end{equation}
 By Claim \ref{c.range of center Lyapunov exponent}, the inequality~(\ref{equa:more sufficient condition}) and  the density of rational numbers on real line, there exist $m,n$ arbitrarily large satisfying the inequality~(\ref{equa:sufficient condition}).
\medskip

 By an arbitrarily $C^1$ small perturbation, the eigenvalues of the periodic orbit $P_1$ are  of multiplicity one (might have complex eigenvalue). Since the period of $P_1$ can be chosen arbitrarily large, by  Lemma~\ref{l.get simple spectrum}, after another small Franks-type  perturbation,  we have that the periodic orbit $P_1$ has simple spectrum.

 %One can see from the proof of \cite[Lemma 4.16]{BDP}, for any $\delta>0$ and $\rho\in(0,1)$,
 %there exists a hyperbolic periodic orbit $\tilde{P_1}$ such that
 %\begin{itemize}
 %\item the periodic orbit $\tilde{P_1}$ is $(\delta,\rho)$-good approximation of $P_1$;
 %\item For each $i=1,\cdots, n$, we have that
 %$$|\chi_{i}(\tilde{P_1})-\chi_i(P_1)|<\delta.$$
 %\end{itemize}
% Hence, when we choose $\delta$ small enough and $\rho$ close to $1$ enough, the periodic orbit $\tilde{P_1}$
%satisfies the conclusion of Lemma~\ref{main lemma}.

 %Since $m,n$ can be chosen arbitrarily large,  $\pi(P_1)$ can be arbitrarily large. Hence, by Franks Lemma,
 %there exists an arbitrarily $C^1$ small Franks type perturbation of $g$ (we still denote it as $g$) supported on
 %an arbitrarily small neighborhood of $P_1$ such that $P_1$ has simple spectrum. Then the items $1,3,4$ are satisfied.
%The items $2,5$ are satisfied by the construction of $g$.

This ends the proof of Lemma~\ref{main lemma}.
\endproof
\begin{remark} One can see from the proof of Lemma~\ref{main lemma} that the perturbation is done in very small neighborhood of the heterodimensional cycle $K$.
\end{remark}

Now we are ready to give the proof of Proposition~\ref{proposition}.
\proof[Proof of Proposition~\ref{proposition}]

We can see that the properties stated in Proposition~\ref{proposition} are persistent under $C^1$ small perturbation. Let $\mathcal{R}$ be  the residual subset of $\diff^1(M)$ from Theorem~\ref{Thm:generic properties}. Notice that for any $f\in\mathcal{R}$, by the item~\ref{generic:periodic orbits with simple spectrum} of Theorem~\ref{Thm:generic properties}, there is a periodic orbit $Q_0$ with simple spectrum which is homoclinically related to $Q$. We take $\chi=-\chi_{i+k}(Q_0,f)>0$.

We only need to show that given $f\in\mathcal{R}$, for any $\zeta>0$ and $\gamma>0$, there are a diffeomorphism $g$ which is $\zeta$-$C^1$-close to $f$ and a hyperbolic periodic orbit $P_1$ of $g$, such that the following properties are satisfied:
\begin{enumerate}
\item $g$ coincides with $f$ on $P_0\cup Q_0$;
\item $P_1$ is robustly in the chain recurrence class of $P_g$;
\item $P_1$ has simple spectrum and the Hausdorff distance $\ud_H(P_1,H(P_g,g))<\gamma$;
\item $\chi_{i+k}(P_1,g)<\frac{3}{4}\cdot\chi_{i+k}(P_{0},g)$;
\item $P_1$ is a $(\gamma, 1-\frac{\chi_{i+k}(P_0,f)}{\chi+\chi_{i+k}(P_0,f)})$-good approximation of $P_0$.
\end{enumerate}
Then Proposition~\ref{proposition} can be proved by a standard Baire argument.

By item \ref{generic:robust in a chain class} of Theorem \ref{Thm:generic properties}, we can require that  $\zeta$ is chosen small enough such that after any $\zeta$-perturbation, the continuations of $P$, $Q$, $P_0$,  and $Q_0$ are still robustly in the same chain recurrence class. We take $0<\epsilon<\frac{\zeta}{4}$, then there exist  $T>0$ and $l_0$ satisfying Lemma \ref{lemma:bb}.

\paragraph{Perturb to get a heterodimensional cycle}
Since $H(P,f)$ admits no dominated splitting of index $j$  for any $j\in\{i+1,\cdots,i+k-1\}$,  there is a number $\delta_0\in(0,\frac{\gamma}{10})$  such that  for any compact invariant subset $\Lambda$ of $H(P,f)$, if $\ud_H(\Lambda,H(P,f))<\delta_0$,  then $\Lambda$ admits no $T$-dominated splitting of index $j$  for any $j\in\{i+1,\cdots,i+k-1\}$.

We fix a positive number $\delta<\min\{\delta_0,\frac{1}{4}\chi_{i+k}(P_0,f)\}$ small enough such that the following is satisfied:
\begin{equation}\label{equa:delta}
\frac{\chi_{i+k}(P_0£¬f)+\delta}{-2\chi_{i+k}(Q_0,f)+\chi_{i+k}(P_0,f)-\delta}
<\frac{\chi_{i+k}(P_0,f)}{-\frac{3}{2}\chi_{i+k}(Q_0,f)+\chi_{i+k}(P_0,f)}.
\end{equation}

We take a number $\kappa$ such that $$\kappa\in\Big(\frac{2\chi_{i+k}(P_0,f)-3\chi_{i+k}(Q_0,f)}{3\chi_{i+k}(P_0,f)-3\chi_{i+k}(Q_0,f)}, 1\Big).$$
 We apply the item~\ref{generic:shadowing periodic orbits} of Theorem \ref{Thm:generic properties} to  the constants $\delta$ and $\kappa$, then  there exist two hyperbolic periodic orbits     $P'=\orb(p')$ and $Q'=\orb(q')$  such that:
\begin{itemize}
\item $P^{\prime}$ and
$Q^{\prime}$  are homoclinically related to $P_0$ and $Q_0$ respectively;
 \item Both $P^{\prime}$ and $Q^{\prime}$ are $\delta/2$ dense in $H(P,f)$ and have simple spectrum.

  \item  $P^{\prime}$ is a $(\frac{\gamma}{10}, \kappa)$-good approximation of $P_0$ and $Q^{\prime}$ is $(\frac{\gamma}{10}, \kappa)$-good approximation of $Q_0$ .

 \item For each $j\in\{1,\cdots,d\}$, we have that
 \begin{equation}\label{equa:approxiamted Lyapunov exponent} |\chi_{j}(P^{\prime},f)-\chi_{j}(P_0,f)|<\delta \textrm{ and } | \chi_{j}(Q^{\prime},f)-\chi_{j}(Q_0,f)|<\delta.
 \end{equation}
\item Both of the periods of $P'$ and $Q'$ are larger than $l_0$.
\end{itemize}

By item~\ref{generic:get cycle} of Theorem~\ref{Thm:generic properties}, we can do an arbitrarily $C^1$ small perturbation,  keeping $P^{\prime}$ and $Q^{\prime}$ homoclinically related to $P$ and $Q$ respectively and without changing the Lyapunov exponents of $P^{\prime}$ and $Q^{\prime}$,  such that  $P^{\prime}$ and $Q^{\prime}$ form a  partially hyperbolic heterodimensional cycle.  For simplicity, we still denote this diffeomorphism as $f$.

Notice that the periodic orbits $P^{\prime}$ and $Q^{\prime}$ have no $T$-domination of index $j$ for any $j\in\{i+1,\cdots,i+k-1\}$.

\paragraph{Equalize the center Lyapunov exponents of both $P'$ and $Q'$}

By Lemma~\ref{lemma:bb} and Remark \ref{r.BB},
there exist   $\pi(P^{\prime})$ one-parameter families $\{(A_{l,t})_{t\in[0,1]}\}_{l=0}^{\pi(P^{\prime})-1}$
and $\pi(Q^{\prime})$ one-parameter families $\{(B_{m,t})_{t\in[0,1]}\}_{m=0}^{\pi(Q^{\prime})-1}$ in $GL(d,\mathbb{R})$
%for any integer $l\in[0,\pi(P^{\prime})-1]$ and any integer $m\in[0,\pi(Q^{\prime})-1]$,
such that:

\begin{itemize}\item $A_{l,0}={\rm Df}(f^{l}(p^{\prime}))$ and $B_{m,0}={\rm Df}(f^{m}(q^{\prime}))$, for any $l,m$;
\item $\norm{A_{l,t}-{\rm Df}(f^{l}(p^{\prime}))}<\epsilon$ and  $\norm{A_{l,t}^{-1}-{\rm Df}^{-1}(f^{l+1}(p^{\prime}))}<\epsilon$, for any $t\in[0,1]$;
\item  $\norm{B_{m,t}-{\rm Df}(f^{m}(q^{\prime}))}<\epsilon$ and  $\norm{B_{m,t}^{-1}-{\rm Df}^{-1}(f^{m+1}(q^{\prime}))}<\epsilon$, for any $t\in[0,1]$;
\item $A_{\pi(P^{\prime})-1,t}\comp\cdots \comp A_{0,t}$ and $B_{\pi(Q^{\prime})-1,t}\comp\cdots \comp B_{0,t}$ are hyperbolic, for any $t\in[0,1]$;
\item   For any integer $s\in [1,i]\cup[i+k+1,d]$, we have that
       $$\chi_{s}(A_{\pi(P^{\prime})-1,t}\comp\cdots \comp A_{0,t})=\chi_{s}(P^{\prime},f) \textrm{ and  }\chi_{s}(B_{\pi(Q^{\prime})-1,t}\comp\cdots \comp B_{0,t})=\chi_{s}(Q^{\prime},f);$$
\item $\chi_{i+1}(A_{\pi(P^{\prime})-1,1}\comp\cdots\comp A_{0,1})=\chi_{i+k}(A_{\pi(P^{\prime})-1,1}\comp\cdots \comp A_{0,1})=\frac{1}{k}\sum_{j=i+1}^{i+k}\chi_{j}(P^{\prime},f)$
    \item $\chi_{i+1}(B_{\pi(Q^{\prime})-1,1}\comp\cdots \comp B_{0,1})=\chi_{i+k}(B_{\pi(Q^{\prime})-1,1}\comp\cdots\comp B_{0,1})=\frac{1}{k}\sum_{j=i+1}^{i+k}\chi_{j}(Q^{\prime},f).$
\end{itemize}

Fix a small number $\eta>0$.
Since $P^{\prime}$ and $Q^{\prime}$ form a heterodimensional cycle and $P^{\prime}$ is homoclinically related to $P_0$,
 there exist four points $x, y, z, w\in M$ such that
\begin{itemize}
\item $$x\in W^{s}_{\eta}(P^{\prime})\cap W^{u}(Q^{\prime}) \textrm{ and } y\in W^{u}_{\eta}(P^{\prime})\cap W^{s}(Q^{\prime});$$
 \item $$z\in W^{s}_{\eta}(P^{\prime})\cap W^{u}(P_0) \textrm{ and } w\in W^{u}_{\eta}(P^{\prime})\cap W^{s}(P_0).$$
\end{itemize}
We take $K^s=\{x,z\}$ and $K^{u}=\{y,w\}$, and we choose a small neighborhood $U$ of $P^{\prime}$ such that $U$ is disjoint from $\orb^{-}(x)\cup \orb^{-}(z)$, $\orb^{+}(y)\cup \orb^{+}(w)$, $Q^{\prime}$, $Q_0$ and two homoclinic orbits between $Q^{\prime}$ and $Q_0$, whose $\omega$-limit sets are $Q'$ and $Q_0$ respectively. By Lemma~\ref{lemma:gourmelon}, there exists an $\epsilon$ perturbation $g_1$ whose support is contained in $U$ such that
 \begin{itemize}
 \item $g_1$ keeps $P^{\prime}$;
  \item ${\rm Dg}_1(f^{l}(p^{\prime}))=A_{l,1}$, for any $l=0,\cdots,\pi(P^{\prime})-1$;
\item $x\in W^{s}(P^{\prime}, g_1)\cap W^{u}(Q^{\prime},g_1)$ and $y\in W^{u}(P^{\prime},g_1)\cap W^{s}(Q^{\prime},g_1)$.
\item $W^{s} (P^{\prime},g_1)$ intersects $ W^{u}(P_0,g_1)$ transversely at the point $z$
         and $W^{u}(P^{\prime},g_1)$ intersects  $W^{s}(P_0,g_1)$ transversely at the point $w$.
\end{itemize}
Following the same  way above, we choose a small neighborhood $V$ of $Q^{\prime}$ which is disjoint from certain homoclinic intersection and some orbit segments, and we apply Lemma \ref{lemma:gourmelon}.
At the end, we get an $\epsilon$ perturbation $g_2$ of $g_1$
such that
\begin{itemize}
\item For diffeomorphism $g_2$, the periodic orbits $P^{\prime}$, $P$, $Q$  and $Q^{\prime}$ are robustly in the same chain recurrence  class;
\item $P^{\prime}$ and $Q^{\prime}$ form a partially hyperbolic heterodimensional cycle;
\item ${\rm Dg}_2(f^{l}(p^{\prime}))=A_{l,1}$ and ${\rm Dg}_2(f^{m}(q^{\prime}))=B_{m,1}$, for any integer $l\in[0,\pi(P^{\prime})-1]$ and $m\in[0,\pi(Q^{\prime})-1]$.
\end{itemize}

To sum up, the diffeomorphism $g_2$   is $2\epsilon$-$C^1$-close to $f$ and satisfies that:
\begin{itemize}
\item[S1.] $g_2$ coincides with $f$ on $P_0\cup Q_0\cup P'\cup Q'$;
\item[S2.] $\chi_j(P_0,g_2)=\chi_j(P_0,f)$ and $\chi_j(Q_0,g_2)=\chi_j(Q_0,f)$, for any $j=1,2,\cdots,d$;
\item[S3.] $\chi_j(P',g_2)=\chi_j(P',f)$ and $\chi_j(Q',g_2)=\chi_j(Q',f)$, for any $j\in [1,i]\cup [i+k+1,d]$;
\item[S4.] $\chi_{i+1}(P',g_2)=\chi_{i+k}(P',g_2)$ and $\chi_{i+1}(Q',g_2)=\chi_{i+k}(Q',g_2)$;
\item[S5.] $\sum_{j=i+1}^{i+k}\chi_{j}(P^{\prime},g_2)=\sum_{j=i+1}^{i+k}\chi_{j}(P^{\prime},f)$  and $\sum_{j=i+1}^{i+k}\chi_{j}(Q^{\prime},g_2)=\sum_{j=i+1}^{i+k}\chi_{j}(Q^{\prime},f)$.
\end{itemize}

As a consequence, the diffeomorphism $g_2$ satisfies the assumptions of Lemma~\ref{main lemma}.
\paragraph{Construction of the periodic orbit $P_1$}
By Lemma \ref{main lemma}, there exist a diffeomorphism  $g$ which is  $\epsilon$-$C^1$-close to $g_2$, hence is $\zeta$-$C^1$-close to $f$, and a hyperbolic periodic orbit $P_1$ of index $i$ for the diffeomorphism $g$ such that
 \begin{itemize}
 \item $\chi_{i+k}(P_1,g)<\frac{1}{2}\cdot\chi_{i+k}(P^{\prime},g)$;
 \item  $P_1$ has  simple spectrum;
 \item $g$ coincides with $g_2$ in a small neighborhood of $P_0\cup Q_0$;
 \item $g$ and ${\rm Dg}$ coincide  with $g_2$ and ${\rm Dg}_2$ on $P^{\prime}\cup Q^{\prime}$ respectively;
 \item $P_1$ is a $(\frac{\gamma}{10}, 1+\frac{\chi_{i+k}(P^{\prime}, g_2)}{2\chi_{i+k}(Q^{\prime},g_2)-\chi_{i+k}(P^{\prime},g_2)})$-good approximation of $P^{\prime}$.
  \end{itemize}
Moreover, by Remark \ref{r.robustly in the chain class}, we have that $P_1$ is robustly in the same chain class with $P_g$ and $Q_g$.
By the choice of $\delta$ and $\gamma$, we have that $\ud_H(P_1,H(P,g))<\gamma$. Then the items $1,2,3$ are satisfied.

By the properties S4 and S5, we have that
$$0<\chi_{i+k}(P',g_2)<\chi_{i+k}(P',f) \textrm{ and } \chi_{i+k}(Q^{\prime},g_2)<\chi_{i+k}(Q^{\prime},f)<0,$$
 which implies that
 \begin{equation}\label{equa:g leq f}
 \frac{\chi_{i+k}(P^{\prime},g_2)}{-2\chi_{i+k}(Q^{\prime},g_2)+\chi_{i+k}(P^{\prime},g_2)} <\frac{\chi_{i+k}(P',f)}{-2\chi_{i+k}(Q',f)+\chi_{i+k}(P',f)}.
 \end{equation}

By the inequalities~(\ref{equa:delta}), ~(\ref{equa:approxiamted Lyapunov exponent}) and ~(\ref{equa:g leq f}), we have that
 \begin{align*}
  \frac{\chi_{i+k}(P^{\prime},g_2)}{-2\chi_{i+k}(Q^{\prime},g_2)+\chi_{i+k}(P^{\prime},g_2)}
  &<\frac{\chi_{i+k}(P_0,f)}{-\frac{3}{2}\chi_{i+k}(Q_0,f)+\chi_{i+k}(P_0,f)}.
 \end{align*}
Recall that  $P^{\prime}$ is a $(\frac{\gamma}{10}, \kappa)$ good approximation of $P_0$ (for the diffeomorphisms $f$ and $g$),
 hence  we have that $P_1$ is a
 $\big(\gamma, \kappa\cdot\big(1+\frac{\chi_{i+k}(P_0,f)}{\frac{3}{2}\chi_{i+k}(Q_0,f)-\chi_{i+k}(P_0,f)}\big)\big)$-
 good approximation of $P_0$.

By the choice of $\kappa$, we have that

\begin{align*}
  &\kappa\cdot\Big(1+\frac{\chi_{i+k}(P_0,f)}{\frac{3}{2}\chi_{i+k}(Q_0,f)-\chi_{i+k}(P_0,f)}\Big)
  \\
  &>\Big(\frac{\chi_{i+k}(P_0,f)-\frac{3}{2}\chi_{i+k}(Q_0,f)}{\frac{3}{2}\chi_{i+k}(P_0,f)-\frac{3}{2}\chi_{i+k}(Q_0,f)}\Big)
  \cdot\Big(1+\frac{\chi_{i+k}(P_0,f)}{\frac{3}{2}\chi_{i+k}(Q_0,f)-\chi_{i+k}(P_0,f)}\Big)
  \\
  &=\frac{\chi_{i+k}(Q_0,f)}{\chi_{i+k}(Q_0,f)-\chi_{i+k}(P_0,f)}
  \\
  &=1+\frac{\chi_{i+k}(P_0,f)}{\chi_{i+k}(Q_0,f)-\chi_{i+k}(P_0,f)}
  \\
  &=1-\frac{\chi_{i+k}(P_0,f)}{\chi+\chi_{i+k}(P_0,f)}
\end{align*}
Hence, $P_1$ is a $\big(\gamma, 1-\frac{\chi_{i+k}(P_0,f)}{\chi+\chi_{i+k}(P_0,f)}\big)$-good approximation of $P_0$. This implies that the item $5$ is satisfied.

  Besides, by the choice of $\delta$, we have the following estimation for the maximal center Lyapunov exponent of $P_1$:
  \begin{displaymath} \chi_{i+k}(P_1,g)<\frac{1}{2}\cdot\chi_{i+k}(P^{\prime},g)=\frac{1}{2}\cdot\frac{1}{k}\cdot\sum_{j=1}^{k}\chi_{i+j}(P^{\prime},f)
  <\frac{1}{2}\cdot\frac{1}{k}\cdot\sum_{j=1}^{k}\chi_{i+j}(P_0,g)+\delta
  \leq\frac{3}{4}\cdot
  \chi_{i+k}(P_0,g).
  \end{displaymath}
  Hence the item $4$ is satisfied. This ends the proof of Proposition~\ref{proposition}.
\endproof

%%%%%%%%%%%%%%%%%%%%%%%%%%%%%%%%%%%%%%%%%%%%%%%%%%%%%%%%%%%%%%%%%%%%%%%%%%%%%%%%
%%%%%%%%%%%%%%%%%%%%%%%%%%%%%%%%%%%%%%%%%%%%%%%%%%%%%%%%%%%%%%%%%%%%%%%%%%%%%%%%
\section{Ergodic measure with multi-zero   Lyapunov exponents for the case  controlled by Jacobian: Proof of Theorem~\ref{Thm: B}}\label{Section:proof of thm B}
%%%%%%%%%%%%%%%%%%%%%%%%%%%%%%%%%%%%%%%%%%%%%%%%%%%%%%%%%%%%%%%%%%%%%%%%%%%%%%%%%
%%%%%%%%%%%%%%%%%%%%%%%%%%%%%%%%%%%%%%%%%%%%%%%%%%%%%%%%%%%%%%%%%%%%%%%%%%%%%%%%%

Consider a diffeomorphism $f\in\diff^1(M)$  and a homoclinic class  $H(P,f)$ admitting a dominated splitting of the form $T_{H(P,f)}M=E\oplus F\oplus G$. We denote by $k=\dim(F)$.  For any periodic orbit $Q=\orb(q)$ contained in $H(P,f)$, the \emph{mean Lyapunov exponent along the bundle $F$} of $Q$ is defined as
$$L^F(Q,f)=\frac{1}{k\cdot\pi(Q)}\log{|\Jac({\rm Df}^{\pi(Q)}|_{F(q)})|}.$$
Notice that $L^F(Q,f)$ is the average of the Lyapunov exponents of $Q$ along the bundle $F$.

%%%%%%%%%%%%%%%%%%%%%%%%%%%%%%%%%%%%%%%%%%%%%%%%%%%%%%%%%%%%%%%%%%%%%%%%%%%%%%%
\subsection{ Proof of Theorem~\ref{Thm: B}}
%%%%%%%%%%%%%%%%%%%%%%%%%%%%%%%%%%%%%%%%%%%%%%%%%%%%%%%%%%%%%%%%%%%%%%%%%%%%%%%
The main ingredient for the proof of Theorem~\ref{Thm: B} is the following proposition.

\begin{proposition}\label{p.jacobi} For generic diffeomorphism $f\in\diff^1(M)$, consider a hyperbolic periodic orbit $P$ of index $i$. Assume the homoclinic class $H(P,f)$ admits  a dominated splitting $T_{H(P,f)}M=E\oplus F\oplus G$, such that $\dim(E)=i$. Assume, in addition, that we have the following:
\begin{itemize}
\item[--] $H(P,f)$ contains a hyperbolic periodic orbit $Q=\orb(q)$, whose index is no larger than $\dim(E\oplus F)$, such that $$|\Jac({\rm Df}^{\pi(Q)}|_{F(q)})|<1;$$
\item[--] the center bundle $F$ has no finer dominated splitting.
\end{itemize}
Then there exists a constant  $\rho\in(0,1)$ which only depends on $Q$,
  such that  for any hyperbolic periodic orbit $P_{0}$ with simple spectrum, which is  homoclinically related to $P$,  and  any  $\gamma>0$,
there exists a hyperbolic periodic point $P_1$ with simple spectrum such that:
\begin{enumerate}
\item $L^F(P_1,f)<\rho \cdot L^{F}(P_0,f)$;
\item $P_1$ is homoclinically related to $P$ and  is $\gamma$ dense inside $H(P,f)$;
\item $P_1$ is $\big(\gamma, 1-\frac{L^F(P_0,f)}{L^F(P_0,f)-L^F(Q,f)}\big)$ good approximation of $P_{0}$.
\end{enumerate}
\end{proposition}

The proof of Proposition~\ref{p.jacobi} is left to the next subsection. Now, we follow the strategy of the Proof of Theorem \ref{Thm: A} to  give the proof of Theorem \ref{Thm: B}.

\proof[Proof of Theorem \ref{Thm: B}]
We denote by $i=\dim(E)$.
By  item \ref{generic:homoclinic class coinsides},  item \ref{generic:periodic orbits with simple spectrum} and item \ref{generic:shadowing periodic orbits}  of Theorem \ref{Thm:generic properties},
we can assume that $P$ is of index $i$  and    has simple spectrum.
 Let $\rho\in(0,1)$ be the number in Proposition \ref{p.jacobi}, which only depends on $Q$.

We will inductively get a sequence  of periodic orbits $\{P_n\}$, a sequence of positive numbers $\{\epsilon_n\}$ and a sequence of integers $\{N_n\}$ satisfying the following properties:
 \begin{itemize}    \item $\epsilon_{n}<\frac{1}{2}\epsilon_{n-1}$;
  \item $L^F(P_{n+1})<\rho\cdot L^F(P_n)$;
 \item $P_{n+1}$ is homoclinically related to $P$ and is $\epsilon_n$ dense inside $H(p,f)$;

 \item $P_{n+1}$ is $(\epsilon_{n},1-\frac{2L^F(P_n,f)}{2L^F(P_n,f)-L^F(Q,f)})$ good approximation of $P_n$;
     \item For any point $x\in B_{2\epsilon_{n}}(P_n)\cap H(P,f)$, we have that 
      $$0<\frac{1}{N_{n}}\log \um({\rm Df}^{N_{n}}|_{F(x)})\leq \frac{1}{N_{n}}\log\norm{{\rm Df}^{N_{n}}|_{F(x)}}<2\chi_{i+k}(P_n).$$
     
 \end{itemize}

\paragraph{Choice of $P_0$, $N_0$ and $\epsilon_0$}
 Let $ P_0=P $, then there exists  an integer $N_{0}$ large enough such that
for any $y\in P_{0}$, we have that
 $$0<\frac{1}{N_0}\log\um({\rm Df}^{N_{0}}|_{F(y)})\leq \frac{1}{N_{0}}\log\norm{{\rm Df}^{N_{0}}|_{F(y)}}<\frac{3}{2}\chi_{i+k}(P_0) .$$

By the uniform continuity of the functions $\log\norm{{\rm Df}^{N_{0}}|_{F(x)}}$ and $ \log\um({\rm Df}^{N_{0}}|_{F(x)})$,
there exists a number $\epsilon_{0}>0$ such that for any point $x\in B_{2\epsilon_{0}}(P_0)\cap H(P,f)$, we have that
  $$0<\frac{1}{N_{0}}\log\um({\rm Df}^{N_{0}}|_{F(x)})\leq\frac{1}{N_{0}}\log\norm{{\rm Df}^{N_{0}}|_{F(x)}}< 2\chi_{i+k}(P_0).$$

\paragraph{Construct $P_n$, $N_n$ and $\epsilon_n$ inductively}
     Assume that $P_i$, $N_i$ and $\epsilon_i$ are already defined for any $i\leq n$.
      We apply $P_n$ and $\epsilon_n$  to the Proposition \ref{p.jacobi},  then we get a periodic orbit $P_{n+1}$  which is homoclinically related to $P_n$  such that
      \begin{itemize}\item  $L^F(P_{n+1})<\rho\cdot L^F(P_n)$;
      \item $P_{n+1}$ is $\epsilon_n$ dense in $H(p,f)$;
       \item $P_n$ is $(\epsilon_{n}, 1-\frac{2L^F(P_n,f)}{2L^F(P_n,f)-L^F(Q,f)})$ good approximation of $P_n$.
      \end{itemize}
Then there exists an integer  $N_{n+1}$ large enough such that
for any $y\in P_{n+1}$, we have that
 $$0<\frac{1}{N_{n+1}}\log\um({\rm Df}^{N_{n+1}}|_{F(y)})\leq \frac{1}{N_{n+1}}\log\norm{{\rm Df}^{N_{n+1}}|_{F(y)}}<\frac{3}{2}\chi_{i+k}(P_{n+1}).$$
By   the uniform continuity of the functions $\log\norm{{\rm Df}^{N_{n+1}}|_{F(x)}}$ and $ \log\um({\rm Df}^{N_{n+1}}|_{F(x)})$,
 there exists a number $\epsilon_{n+1}\in(0,\frac{1}{2}\epsilon_n]$ such that  
 for any point $x\in B_{2\epsilon_{n+1}}(P_{n+1})\cap H(P,f)$, we have
 \begin{displaymath} 0<\frac{1}{N_{n+1}}\cdot\log \um({\rm Df}^{N_{n+1}}|_{F(x)})\leq\frac{1}{N_{n+1}}\log\norm{{\rm Df}^{N_{n+1}}|_{F(x)}}< 2\chi_{i+k}(P_{n+1}).
   \end{displaymath}

\paragraph{End of proof of Theorem~\ref{Thm: B}}
Since  $1-\frac{2L^F(P_n,f)}{2L^F(P_n,f)-L^F(Q,f)}$ exponentially tends to $1$ and $\sum_n\epsilon_n$ converges,
by Lemma~\ref{limit}, the sequence of ergodic measures $\delta_{P_n}$ converges  to an ergodic measure $\nu$ whose support is $H(p,f)$.

 \begin{claim}\label{c.zero Lyapunov exponent jacobi}
 The Lyapunov exponents of $\nu$ along the center bundle $F$ are all zero.
 \end{claim}
Notice that $\chi_{i+k}(P_n)\leq k\cdot L^F(P_n)\leq k\cdot \rho^{n} L^F(P_0)$.  The  proof of Claim~\ref{c.zero Lyapunov exponent jacobi} follows the proof of the Claim~\ref{c.zero Lyapunov exponent}. The only difference is that we control the sum of the center Lyapunov exponents by the function $\frac{1}{N_{n}}\log\norm{{\rm Df}^{N_{n }}|_{F}}$ instead of  the function $L_{d-i}^{N_n}-L_{d-i-k}^{N_n}$.
\medskip

This ends the proof of Theorem~\ref{Thm: B}.

\endproof
Now it remains to prove Proposition~\ref{p.jacobi}
%%%%%%%%%%%%%%%%%%%%%%%%%%%%%%%%%%%%%%%%%%%%%%%%%%%%%%%%%%%%%%%%%%%%%%%%%%%%%%
\subsection{Good approximation with weaker center Jacobian: Proof of Proposition~\ref{p.jacobi}}
%%%%%%%%%%%%%%%%%%%%%%%%%%%%%%%%%%%%%%%%%%%%%%%%%%%%%%%%%%%%%%%%%%%%%%%%%%%%%%

The proof of Proposition~\ref{p.jacobi} is based on the following perturbation lemma:
 \begin{lemma}\label{l.jacobi}
 Let  $P$ and $Q$ be two hyperbolic periodic orbits of $f\in \diff^1(M)$ with different indices.
 Assume that
 \begin{itemize}
 \item[--] $Q$ and $P$ form a  partially hyperbolic heterodimensional cycle $K$.
       In other words, $K$ admits a partially hyperbolic splitting of the form
        $$T_{K}M=E^s\oplus E^c\oplus E^u,$$
         where $\dim(E^s)=\ind(P)$ and $\dim(E^s\oplus E^c)=\ind(Q)$;
 \item[--]  there exists another dominated splitting over $K$ of the form $$T_{K}M=E^s\oplus F\oplus G$$
           such that   $\dim(F)\geq\dim(E^c)$;
 \item[--] all the Lyapunov exponents of $Q$ along $E^c$ are equal.
  \item[--] all the Lyapunov exponents of $P$ along $E^c$ are equal and are larger than  $L^F(P,f)/2$;
 \item[--] $L^{F}(Q,f)<0.$
 \end{itemize}
Then there exists a number $\rho\in(0,1)$ which only depends on $Q$,    such that
 for any $\gamma>0$  and   any $C^1$ neighborhood $\mathcal{U}$ of $f$, there exists $g\in\mathcal{U}$ together with a hyperbolic periodic orbit $P^{\prime}$ of index $\ind(P)$, with simple spectrum such that
 \begin{enumerate}
 \item $g=f$ and ${\rm Dg}={\rm Df}$ on $P\cup Q$ ;
 \item $L^F(P^{\prime},g)
     <\rho \cdot L^F(P,g)$;
 \item  $P^{\prime}$ is $\big(\gamma, 1-\frac{ L^F(P,g)}{L^F(P,g)-L^F(Q,g)})$ good approximation of $P$;
 \item $W^{ss}(P^{\prime},g)$ has transverse intersections with $W^{u}(P,g)$ and $W^{uu}(P^{\prime},g)$ has transverse intersections with $W^s(Q,g)$,  corresponding  to the partially hyperbolic splitting $T_{K_g}M=E^s\oplus E^c\oplus E^u$.
 \end{enumerate}
 \end{lemma}
 \begin{remark}\label{r.jacobian robustly in the same chain class}
 \begin{enumerate}
 \item Once again, if $P$ and $Q$ are robustly in the same chain recurrence class, the fourth item above implies that $P^{\prime}$, $P$ and $Q$ are robustly in the same chain recurrence class;
 \item Actually, the constant $\rho$ is only and continuously depends on the  mean Lyapunov exponent of $Q$ along the bundle $E^c$ and the mean  Lyapunov exponent of $Q$ along $F$.
 \end{enumerate}
 \end{remark}
 The idea of the proof of Lemma~\ref{l.jacobi} is that we mix two hyperbolic periodic orbits with different sign of mean Lyapunov exponents to get a new hyperbolic periodic orbit with weaker mean  Lyapunov exponent along the bundle $F$.

Similar to Section~\ref{Section:proof of thm A},  we complete the proof of Proposition~\ref{p.jacobi} by proving Lemma~\ref{l.jacobi}.
To prove Lemma~\ref{l.jacobi}, we  first follow the strategy of the proof of  Lemma \ref{main lemma} to linearize the system in a small neighborhood of the cycle $K$ by an arbitrarily small perturbation,  then by another arbitrarily small perturbation, we get a periodic orbit.
  At the end, we will adjust the time of periodic orbit staying close to $P$ and $Q$ respectively.

 \proof[Proof of Lemma~\ref{l.jacobi}]
 By the assumption, we can denote by $\log\mu$ and $\log\lambda$ the Lyapunov exponents of $P$ and $Q$ along $E^c$ respectively.
 Then we have that
 $$\frac{L^F(P,f)}{2}< \log\mu< L^F(P,f) \textrm{ and } \log\lambda<L^F(Q,f).$$
 Denote by $P=\orb(p,f)$ and $Q=\orb(q,f)$.

We fix a small number $\gamma>0$ and a neighborhood $\mathcal{U}$ of $f$.
 There exists $\epsilon>0$ such that the $\epsilon$ neighborhood of $f$ is contained in $\mathcal{U}$.
There is a small number $0<\theta<1$, such that  for any $h\in\mathcal{U}$, if $d(z_1,z_2)<\theta\cdot\gamma$,
we have that
$$d(h^i(z_1),h^i(z_2))<\frac{\gamma}{2}, \textrm{  for any $i\in [-\pi(P),\pi(P)]$}.$$
We take two neighborhoods $U_P$ and $U_Q$ of $P$ and $Q$ respectively,
such that $U_P$ is contained in the $\theta\cdot\gamma$-neighborhood of $P$ and is disjoint from $U_Q$.

\paragraph{Construction of the periodic orbit $P_1$}
 Similar  to the proof of Lemma \ref{main lemma}, consider  the splitting $T_KM=E^s\oplus E^c\oplus E^{u}=E^s\oplus F\oplus G$ and the two neighborhoods $U_P$ and $U_Q$, by Theorem~\ref{Thm:transition} and Remark~\ref{r.transition},
 there are two matrices $T_0, T_1$ and two  positive  integers $t_0,t_1$ such that  for any two integers $m$ and $n$,
  there exist $g\in\mathcal{U}$ and a hyperbolic periodic orbit $P_1=\orb(p_1,g)$ satisfying s that:
 \begin{itemize}
 \item[--] $g=f$ and ${\rm Dg}={\rm Df}$ on $P\cup Q$,
 \item[--] The linear maps   $T_0: T_pM \mapsto T_qM$ and $T_1: T_qM\mapsto T_pM$ preserve  the two dominated splittings,
  \item[--] $\pi(P_1)=m\pi(P)+n\pi(Q)+t_0+t_1$,
 \item[--] ${\rm Dg}^{\pi(P_1)}(p_1)$ is conjugate to  $T_1\circ {\rm Df}^{n\pi(Q)}(q) \circ T_0\circ {\rm Df}^{m\pi(P)}(p),$
 \item[--] $\#(P_1\cap U_P)\geq m\pi(P)$ and $\#(P_1\cap U_Q)\geq n\pi(Q)$,
 \item[--] $W^{ss}(P_1,g)\pitchfork W^{u}(P,g)$ and $W^{uu}(P_1,g)\pitchfork W^s(Q,g)$ corresponding to the splitting $T_{K_g}M=E^s\oplus E^c\oplus E^u$.
\end{itemize}

 As a consequence of the first item above, we have that
 \begin{itemize}\item  $L^F(P,f)=L^F(P,g)$ \textrm{ and }$L^F(Q,f)=L^F(Q,g);$
 \item $\chi_j(P,f)=\chi_j(P,g) $  and $\chi_j(Q,f)=\chi_j(Q,g)$,  for any $j=1,\cdots,d$.
  \end{itemize}
For simplicity, we   denote them by  $L^F(P)$, $L^F(Q)$, $\chi_j(P)$ and  $\chi_j(Q)$,.

 Since $T_0$ and $T_1$ preserve the dominated splittings, by choosing the proper coordinates, we assume that, corresponding to the two splittings, the two matrices $T_1$ and $T_0$ have the following forms respectively:
 \begin{displaymath}
       \mathbf{T_0}=
        \left(\begin{array}{ccc}
        D_{s}&0&0\\
        0&D_{c}&0\\
        0&0&D_{u}
        \end{array}\right)=
        \left(\begin{array}{ccc}
        D_{s}&0&0\\
        0&D_{F}&0\\
        0&0&D_{G}
        \end{array}\right),
     \end{displaymath}
    \begin{displaymath}
       \mathbf{T_1}=
        \left(\begin{array}{ccc}
        C_{s}&0&0\\
        0&C_{c}&0\\
        0&0&C_{u}
        \end{array}\right)=
        \left(\begin{array}{ccc}
        C_{s}&0&0\\
        0&C_{F}&0\\
        0&0&C_{G}
        \end{array}\right).
    \end{displaymath}

 Then we have that
 ${\rm Dg}^{\pi(P_1)}(p_1)|_{E^c}$ is conjugate to $$C_c\comp {\rm Df}^{n\pi(Q)}(q)|_{E^c}\comp D_{c} \comp {\rm Df}^{m\pi(P)}(p)|_{E^c},$$
 and
 ${\rm Dg}^{\pi(P_1)}(p_1)|_F$ is conjugate to $$C_F\comp {\rm Df}^{n\pi(Q)}(q)|_{F}\comp D_{F} \comp {\rm Df}^{m\pi(P)}(p)|_{F}.$$

 \paragraph{Choice of $m$, $n$ and $\rho$}
 We will adjust $m$ and $n$ to get a periodic orbit satisfying the conclusion of Lemma \ref{l.jacobi}. Let $\eta>0$ be a small number which will be decided later.
 \begin{claim}\label{c.jacobian and Lyapunov exponent}There exists an integer $N_{\eta}$  such that for any $m,n\geq N_{\eta}$, we have that
 \begin{itemize}\item all the  Lyapunov exponents of $P_1$  along the bundle $E^c$ would belong to the interval
 $$\big[\frac{m\pi(P)\cdot\log\mu+n\pi(Q)\cdot\log\lambda}{m\pi(P)+n\pi(Q)}-\eta,\,
\frac{m\pi(P)\cdot\log\mu+n\pi(Q)\cdot\log\lambda}{m\pi(P)+n\pi(Q)}+\eta\big].$$
\item the mean Lyapunov exponent of $P_1$ along the bundle $F$ would belong to the interval
$$\big[\frac{m\pi(p)\cdot L^F(P)+n\pi(q)\cdot L^F(Q)}{m\pi(P)+n\pi(Q)}-\eta,\,
\frac{m\pi(P)\cdot L^F(P)+n\pi(Q)\cdot L^F(Q)}{m\pi(P)+n\pi(Q)}+\eta\big].$$
 \end{itemize}
 \end{claim}
 The proof of Claim \ref{c.jacobian and Lyapunov exponent} is just like the proof of Claim \ref{c.range of center Lyapunov exponent} and we omit the proof here.
 \medskip

To guarantee the item 2 and that $P_1$ has the same index as $P$, we only need to require that there exists a number $\rho\in(0,1)$
 which will be decided later, such that:
%\begin{equation}\label{eq:index}
$$\frac{m\pi(P)\cdot L^F(P)+n\pi(Q)\cdot L^F(Q)}{m\pi(P)+n\pi(Q)}+\eta<\rho\cdot L^F(P)
%\end{equation}
\textrm{\; and\; }
%\begin{equation}\label{eq:mean jacobi}
\frac{m\pi(P)\cdot\log\mu+n\pi(Q)\cdot\log\lambda}{m\pi(P)+n\pi(Q)}-\eta>0,$$
%\end{equation}
%The inequality~(\ref{eq:index}) is equivalent to
%\begin{equation}\label{eq:index imply}
which are equivalent to
$$\frac{\eta-\log\lambda}{\log\mu-\eta}<\frac{m\pi(P)}{n\pi(Q)}<\frac{\rho \cdot L^F(P)-L^F(Q)-\eta}{L^{F}(P)-\rho\cdot L^F(P)+\eta}.$$
%$\end{equation}
%The inequality~(\ref{eq:mean jacobi}) is equivalent to
%\begin{equation}\label{eq:mean jacobi imply}
%\frac{m\pi(P)}{n\pi(Q)}>\frac{\eta-\log\lambda}{\log\mu-\eta}.
%\end{equation}
%If the inequality $$\frac{\rho L^F(P)-L^F(Q)}{L^{F}(P)-\rho L^F(P)}>\frac{-\log\lambda}{\log\mu}$$ is satisfied for some $\rho\in(0,1)$,
\begin{claim}\label{c.choice of rho} There exists $\rho\in(0,1)$ such that the following inequality is satisfied:
$$\frac{-\log\lambda}{\log\mu}< \frac{\rho \cdot L^F(P)-L^F(Q)}{L^{F}(P)-\rho\cdot L^F(P)}.$$
\end{claim}
\proof  The proof consists in solving the following inequality:
%\begin{equation}
$$\frac{\rho\cdot L^F(P)-L^F(Q)}{L^{F}(P)-\rho \cdot L^F(P)}>\frac{-\log\lambda}{\log\mu},$$
%\end{equation}
which is equivalent to
%\begin{equation}
$$\rho>\frac{L^F(Q)\cdot \log\mu-\log\lambda \cdot L^F(P)}{L^F(P)\cdot\log\mu-\log\lambda\cdot  L^F(P)}.$$
%\end{equation}
By assumption  that $\log\mu\in\big(\frac{L^F(P)}{2},L^F(P)\big)$,
   we have the estimation:
%\begin{equation}
$$\frac{L^F(Q)\cdot\log\mu-\log\lambda \cdot L^F(P)}{L^F(P)\cdot\log\mu-\log\lambda\cdot  L^F(P)}
=\frac{L^F(Q)\cdot\frac{\log\mu}{L^F(P)}-\log\lambda }{\log\mu-\log\lambda}
< \frac{\frac{L^F(Q)}{2}-\log\lambda }{-\log\lambda}
=1+\frac{L^F(Q)}{-2\log\lambda}.$$
%\end{equation}
We only need to take $$\rho=1+\frac{L^F(Q)}{-2\log\lambda}\in(0,1).$$
Notice that $\rho$ only depends on $L^F(Q)$ and $\log\lambda$.
\endproof
We fix the value of $\rho$ that we get from Claim \ref{c.choice of rho}, then when $\eta$ is chosen small enough, we have that
$$\frac{\eta-\log\lambda}{\log\mu-\eta}<\frac{\rho\cdot L^F(P)-L^F(Q)-\eta}{L^{F}(P)-\rho\cdot L^F(P)+\eta}.$$
 By the density of rational numbers among $\mathbb{R}$ and Claim \ref{c.jacobian and Lyapunov exponent}, there exist $m$ and $n$ arbitrarily large such that $$\frac{m\pi(P)}{n\pi(Q)}\in\Big(\frac{\eta-\log\lambda}{\log\mu-\eta},  \frac{\rho \cdot L^F(P)-L^F(Q)-\eta}{L^{F}(P)-\rho \cdot L^F(P)+\eta}\Big),$$
which implies that
$P_1$ satisfies the properties of the conclusion of Lemma \ref{l.jacobi} except the item 3.

Now, we only need to check that the choice of $m$ and $n$ guarantees  the item 3.
By the fact that
$$\frac{m\pi(P)}{n\pi(Q)}>\frac{-\log\lambda}{\log\mu}, L^F(P)>\log\mu \textrm{ and }L^F(Q)>\log\lambda,$$
we have that
$$\frac{m\pi(P)}{m\pi(P)+n\pi(Q)}>\frac{\frac{-\log\lambda}{\log\mu}}{\frac{-\log\lambda}{\log\mu}+1}
=\frac{-\log\lambda}{\log\mu-\log\lambda}
> \frac{-L^F(Q)}{L^F(P)-L^F(Q)}=1-\frac{ L^F(P)}{L^F(P)-L^F(Q)}.$$
Hence by taking $m$ and $n$ much larger than $t_0+t_1$, we have that $P^{\prime}$ is $\big(\gamma, 1-\frac{ L^F(P)}{L^F(P)-L^F(Q)})$ good approximation of $P$. Just as the part of the proof of Lemma~\ref{main lemma},  $P$ can be chosen with simple spectrum.
This ends the proof of Lemma \ref{l.jacobi}.
 \endproof
\begin{remark}\label{r.rho}From the proof above, one can see that $\rho$ only depends on $L^F(Q)$ and the average of the
Lyapunov exponents of $Q$ along $E^c$.
 \end{remark}

 Now, we can give the proof of Proposition \ref{p.jacobi} whose proof is quite similar to that of Proposition \ref{proposition}.

 \proof[Proof of Proposition \ref{p.jacobi}]
 We denote by $k=\dim(F)$  and assume that $\ind(Q)=i+k_0$, then we have $0<k_0\leq k$.

 We can see that the properties stated in Proposition~\ref{p.jacobi} are persistent under $C^1$ small perturbation. Let  $\mathcal{R}$ be  the residual subset of $\diff^1(M)$ from Theorem~\ref{Thm:generic properties}. We only need to show that given $f\in\mathcal{R}$, there is $\rho\in(0,1)$ such that for any $\zeta>0$ and $\gamma>0$, there exist a diffeomorphism $g$ which is $\zeta$-$C^1$-close to $f$,  and a $g$ hyperbolic periodic orbit $P_1$ of index $i$ satisfying the followings:
\begin{itemize}
\item[$H_1.$] $g$ coincides with $f$ on $P_0\cup Q$;
\item[$H_2.$]$P_1$ is robustly in the chain recurrence class of $P_g$;
\item[$H_3.$] $P_1$ has simple spectrum and $\ud_H(P_1,H(P_g,g))<\gamma$;
\item[$H_4.$] $L^F(P_1,g)<\rho\cdot L^F(P_{0},g)$;
\item[$H_5.$] $P_1$ is a $(\gamma, 1-\frac{L^F(P_0,g)}{L^F(P_0,g)-L^F(Q,g)})$-good approximation of $P_0$.
\end{itemize}
Then Proposition~\ref{p.jacobi} can be proved by a standard Baire argument.

\paragraph{The previous settings}

By item \ref{generic:robust in a chain class} of Theorem \ref{Thm:generic properties}, we can require that  $\zeta$ is chosen small enough such that after any $\zeta$-perturbation, the continuations of $P_0$, $P$ and $Q$ are still robustly in the same chain recurrence class.

We take $0<\epsilon<\frac{\zeta}{4}$, then there exist  $T>0$ and $l_0$ satisfying Lemma \ref{lemma:bb}.

 We denote by $$\log\lambda=\frac{1}{k_0}\sum_{j=i+1}^{i+k_0}\chi_{j}(Q,f)
 \textrm{ and } \rho_{_0}=1+\frac{L^F(Q,f)}{-2\log\lambda}.$$

Since $H(P,f)$ admits no dominated splitting of index $j$  for any $j\in\{i+1,\cdots,i+k-1\}$,
there is a number $\delta_0\in(0,\frac{\gamma}{10})$ such that
 for any compact invariant subset $\Lambda\subset H(P,f)$, if 
  $\ud_H(\Lambda,H(P,f))<\delta_0,$ 
 then $\Lambda$ admits no $T$-dominated splitting of index $j$  for any $j\in\{i+1,\cdots,i+k-1\}$.

Notice that $\log\lambda<0$, $L^F(Q,f)<0$, $0<L^F(P_0,f)$ and $\rho_{_0}\in(0,1)$, hence we have that
\begin{itemize}
\item $$\rho_{_0}=1+\frac{L^F(Q,f)}{-2\log\lambda}<\frac{1+\rho_{_0}}{2};$$
\item $$\frac{L^F(P_0,f)}{L^F(P_0,f)-L^F(Q,f)}\in\Big(0, \frac{2L^F(P_0,f)}{2L^F(P_0,f)-L^F(Q,f)}\Big).$$
\end{itemize}
As a consequence, we can take a  number $\delta\in(0,\frac{1-\rho_{_0}}{2})$ small enough such that:
\begin{itemize}
\item  $$1+\frac{L^F(Q,f)+\delta}{-2\log\lambda+\delta}<\frac{1+\rho_{_0}}{2};$$
\item $$\frac{L^F(P_0,f)+\delta}{L^F(P_0,f)-L^F(Q,f)-2\delta}\in\Big(0, \frac{2L^F(P_0,f)}{2L^F(P_0,f)-L^F(Q,f)}\Big)$$
\end{itemize}

We take a number $\kappa$ such that
 $$\kappa\in\Big(\frac{2L^F(P_0,f)-L^F(Q,f)}{3L^F(P_0,f)-L^F(Q,f)}, 1\Big).$$
 We apply the item~\ref{generic:shadowing periodic orbits} of Theorem \ref{Thm:generic properties} to  the constants $\delta$ and $\kappa$, then  there exist two hyperbolic periodic orbits   $P'=\orb(p')$ and $Q'=\orb(q')$,  with simple spectrum such that:
\begin{itemize}
\item $P^{\prime}$ and
$Q^{\prime}$  are homoclinically related to $P_0$ and $Q$ respectively;
 \item Both $P^{\prime}$ and $Q^{\prime}$ are $\delta/2$ dense in $H(P,f)$;
  \item  $P^{\prime}$ is a $(\frac{\gamma}{10}, \kappa)$-good approximation of $P_0$ and $Q^{\prime}$ is $(\frac{\gamma}{10}, \kappa)$-good approximation of $Q$;
 \item  $ |L^F(P^{\prime},f)-L^F(P_0,f)|<\delta \textrm{ and } | L^F(Q^{\prime},f)-L^F(Q,f)|<\delta;$
  \item $$\big|\sum_{j=i+1}^{i+k_0}\chi_j(Q^{\prime})-\sum_{j=i+1}^{i+k_0}\chi_j(Q)\big|<\delta;$$
\item Both of the periods of $P'$ and $Q'$ are larger than $l_0$.
\end{itemize}

By item~\ref{generic:get cycle} of Theorem~\ref{Thm:generic properties}, we can do an arbitrarily $C^1$ small perturbation, keeping $P^{\prime}$ and $Q^{\prime}$ homoclinically related to $P$ and $Q$ respectively and without changing the Lyapunov exponents of $P^{\prime}$ and $Q^{\prime}$, such that  $P^{\prime}$ and $Q^{\prime}$ form a  partially hyperbolic heterodimensional cycle.

\paragraph{Equalize the center Lyapunov exponents of both $P'$ and $Q'$}
By Lemma~\ref{lemma:bb}, there exist
 $\pi(P^{\prime})$ one-parameter families $\{(A_{l,t})_{t\in[0,1]}\}_{l=0}^{\pi(P^{\prime})-1}$ and
$\pi(Q^{\prime})$ one-parameter families \\$\{(B_{m,t})_{t\in[0,1]}\}_{m=0}^{\pi(Q^{\prime})-1}$
in $GL(d,\mathbb{R})$ such that:

\begin{itemize}
\item $A_{l,0}={\rm Df}(f^{l}(p^{\prime}))$ and $B_{m,0}={\rm Df}(f^{m}(q^{\prime}))$, for any $l,m$;
\item $\norm{A_{l,t}-{\rm Df}(f^{l}(p^{\prime}))}<\epsilon$ and  $\norm{A_{l,t}^{-1}-{\rm Df}^{-1}(f^{l+1}(p^{\prime}))}<\epsilon$, for any
            $t\in[0,1]$;
\item  $\norm{B_{m,t}-{\rm Df}(f^{m}(q^{\prime}))}<\epsilon$ and  $\norm{B_{m,t}^{-1}-{\rm Df}^{-1}(f^{m+1}(q^{\prime}))}<\epsilon$, for any
            $t\in[0,1]$;
\item $A_{\pi(P^{\prime})-1,t}\comp\cdots \comp A_{0,t}$ and $B_{\pi(Q^{\prime})-1,t}\comp\cdots \comp B_{0,t}$ are hyperbolic, for any
       $t\in[0,1]$;
\item   For any integer $s\in [1,i]\cup[i+k+1,d]$, we have that
       $$\chi_{s}(A_{\pi(P^{\prime})-1,t}\comp\cdots \comp A_{0,t})=\chi_{s}(P^{\prime},f) .$$
\item For any integer $s\in[1,i]\cup[i+k_0+1,d]$, we have that
          $$\chi_{s}(B_{\pi(Q^{\prime})-1,t}\comp\cdots \comp B_{0,t})=\chi_{s}(Q^{\prime},f);$$
\item $\chi_{i+1}(A_{\pi(P^{\prime})-1,1}\comp\cdots\comp A_{0,1})=\chi_{i+k_0}(A_{\pi(P^{\prime})-1,1}\comp\cdots \comp A_{0,1})\in(\frac{1}{2}L^{F}(P^{\prime}), L^F(P^{\prime}));$
\item  $\chi_{i+k_0+1}(A_{\pi(P^{\prime})-1,1}\comp\cdots\comp A_{0,1})=\chi_{i+k}(A_{\pi(P^{\prime})-1,1}\comp\cdots \comp A_{0,1})\geq L^{F}(P^{\prime});$
\item $\chi_{i+1}(B_{\pi(Q^{\prime})-1,1}\comp\cdots \comp B_{0,1})=\chi_{i+k_0}(B_{\pi(Q^{\prime})-1,1}\comp\cdots\comp B_{0,1})=\frac{1}{k_0}\sum_{j=i+1}^{i+k_0}\chi_{j}(Q^{\prime},f).$
\end{itemize}

Similar to the proof of Proposition \ref{proposition}, by Franks-Gourmelon Lemma, there exists an $\epsilon$ perturbation $g_1$ of $f$, which preserves the partially hyperbolic heterodimensional cycle formed by $P^{\prime}$ and $Q^{\prime}$, such that
\begin{itemize}
\item[S1.] $g_1$ coincides with $f$ on $P_0\cup Q\cup P'\cup Q'$;
\item[S2.] $\chi_j(P_0,g_1)=\chi_j(P_0,f)$ and $\chi_j(Q,g_1)=\chi_j(Q,f)$, for any $j=1,2,\cdots,d$;
\item[S3.] $\ind(P^{\prime},g_1)=\ind(P^{\prime},f)$ and $\ind(Q^{\prime},g_1)=\ind(Q^{\prime},f)$;
\item[S5.] $\chi_{i+1}(P',g_1)=\chi_{i+k_0}(P',g_1)\in\big(\frac{1}{2}L^F(P^{\prime},f), L^F(P^{\prime},f)\big)$
\item[S6.] $\chi_{i+1}(Q',g_1)=\chi_{i+k_0}(Q',g_1)$;
\item[S7.] $L^F(P^{\prime},g_1)=L^F(P^{\prime},f)$  and $L^F(Q^{\prime},g_1)=L^F(Q^{\prime},f)$.
\end{itemize}

\paragraph{Construction of the periodic orbit  $P_1$}
By Lemma \ref{l.jacobi} and  Remark \ref{r.jacobian robustly in the same chain class}, we have that  there exist $g\in\diff^1(M)$, which is  $\epsilon$-$C^1$-close to $g_1$ and therefore  is $\zeta$-$C^1$-close to $f$, and a $g$-hyperbolic periodic orbit $P_1$ of index $i$ such that
 \begin{itemize}
 \item  $L^F(P_1,g)<\rho_{_0}\cdot L^F(P^{\prime},g_1)$;
 \item  $P_1$ has  simple spectrum;
 \item  $g$ coincides with $g_1$ on $P_0\cup Q\cup P'\cup Q'$;
 \item  the Lyapunov exponents of $P_0, Q, P' \textrm{ and } Q'$ with respect to $g_1$ are equal  to those with respect to $g$, hence are equal
             to those with  respect to $f$;
  \item $P_1$ is a $(\frac{\gamma}{10}, 1-\frac{L^F(P^{\prime}, f)}{L^F(P^{\prime}, f)-L^F(Q^{\prime}, f)})$-good approximation of $P^{\prime}$;
  \item  $P_1$ is robustly in the same chain class with $P_g$ and $Q_g$.
  \end{itemize}

   By the choice of $\delta$ and $\gamma$, we have that $\ud_H(P_1,H(P,g))<\gamma$. Then the items $H_1, H_2, H_3$ are  satisfied.

By the choice of $\delta$, we have the following estimation for the mean center Lyapunov exponent $L^F(P_1)$ of $P_1$:
  \begin{displaymath} L^F(P_1,g)<\rho_{_0}\cdot L^F(P^{\prime},g_1)=\rho_{_0}\cdot L^F(P^{\prime},f)<\rho_{_0}\cdot L^F(P_0,f)+\delta
  <\frac{1+\rho_{_0}}{2}\cdot L^F(P_0,g).
  \end{displaymath}
  We only need to take $\rho=\frac{1+\rho_{_0}}{2}$, hence item $H_4$ is satisfied.

Besides, by the choice of $\kappa$, we have that
\begin{align*}
  &\kappa\cdot\Big(1-\frac{L^F(P^{\prime}, f)}{L^F(P^{\prime}, f)-L^F(Q^{\prime}, f)}\Big)
  \\
  &>\frac{2L^F(P_0,f)-L^F(Q,f)}{3L^F(P_0,f)-L^F(Q,f)}\cdot \frac{-L^F(Q,f)}{2L^F(P_0,f)-L^F(Q,f)}
  \\
  &=1-\frac{3L^F(P_0,f)}{3L^F(P_0,f)-L^F(Q,f)}.
\end{align*}

Since $P^{\prime}$ is $(\frac{\gamma}{10}, \kappa)$ good approximated of $P_0$,  by the inequality above, $P_1$ is a $\big(\gamma, 1-\frac{2L^F(P_0,f)}{2L^F(P_0,f)-L^F(Q,f)}\big)$-good approximation of $P_0$. Then item $H_5$ is satisfied.

This ends the proof of Proposition~\ref{p.jacobi}.
 \endproof

\section{Partially hyperbolic homoclinic classes with volume expanding center bundle}\label{Section:counterexample}

In this section, we give an example showing that  Corollary~\ref{thm.with dominated splitting} may be not true if there is no periodic orbit of index $\dim(E\oplus F)$. We first give some known results about normally hyperbolic submanifolds in Section~\ref{Section:normally hyperbolic} and the example will be given in Section~\ref{Section:example}.

\subsection{Stability of normally hyperbolic compact manifolds}\label{Section:normally hyperbolic}

 Let $f\in\diff^1(M)$. A compact invariant submanifold without boundary $N$ of $M$ is called \emph{normally hyperbolic}, if there exists a partially hyperbolic splitting of the form $T_{N}M=E^s\oplus TN\oplus E^u$.

  We state a simple version of Theorem 4.1 in \cite{HPS} which gives the stability theorem for normally hyperbolic compact submanifold.
 \begin{theorem}\label{stability of normally hyperbolic} Let $f\in\diff^1(M)$ and $N$ be a compact normally hyperbolic submanifold. We denote by $i:N\mapsto M$  the embedding map from $N$ to $M$.

 There exists a $C^1$ small neighborhood $\mathcal{U}$ of $f$ such that for any $g\in\mathcal{U}$,
 there exists a $C^1$ embedding map $i_{g}:N\mapsto M$, such that $N_{g}=i_{g}(N)$ is $g$-normally hyperbolic. Moreover, $i_{g}$ would tend to $i$ in the $C^1$ topology, if $g$ tends to $f$.
 \end{theorem}

 \begin{remark}\label{Rem:conjugation}
 The map $i_{g}^{-1}|_{N_g}\comp g\comp i_{g}$ is $C^1$-conjugate to the restriction of the map $g$ to $N_g$ and is $C^1$ close to $f$ if $g$ is $C^1$ close to $f$.
 \end{remark}

\subsection{An example}\label{Section:example}

Ch. Bonatti \cite{B} (see also Section 6.2 in \cite{BV}) constructs an open set $\mathcal{U}$ of $C^1$ diffeomorphism  on $\mathbb{T}^3$ such that for any $f\in\mathcal{U}$, we have the following:
  \begin{itemize}\item $f$ is robustly transitive;
   \item There exist a periodic orbit  of index one having a complex eigenvalue and a periodic orbit of index two.
\end{itemize}
 By Theorem 2 and Theorem 4 in \cite{BDP}, for the diffeomorphism $f\in\mathcal{U}$, there exists a partially hyperbolic splitting of the form $T\mathbb{T}^{3}=E^{ss}\oplus E^c$, where $\dim(E^{ss})=1$ and the center bundle $E^c$ is volume expanding without any finer dominated splitting.

 Now, we consider a north-south diffeomorphism $h$ on $\mathbb{S}^1$ such that the expanding rate of $h$ at the source $Q$ is strictly  larger than the norm of $f$. We denote by
 $$\tilde{f}=f\times h: \mathbb{T}^3\times \mathbb{S}^1\mapsto \mathbb{T}^3\times \mathbb{S}^1.$$
 By Theorem~\ref{stability of normally hyperbolic} and continuity of partial hyperbolicity, there exists a $C^1$ neighborhood  $\mathcal{V}$  of $\tilde{f}$ such that  any $\tilde{g}\in\mathcal{V}$ has a partially hyperbolic repelling  set $\Lambda_{\tilde{g}}$ diffeomorphic to $\mathbb{T}^3\times\{Q\}$ admitting a splitting of the form $T_{\Lambda_{\tilde{g}}}\mathbb{T}^4=E^{ss} \oplus E^c\oplus E^{uu}$ where $E^{ss}\oplus E^c=T\Lambda_{\tilde{g}}$. Then by Remark~\ref{Rem:conjugation}, we have that $\tilde{g}|_{\Lambda_{\tilde{g}}}$ is transitive and the dynamics $\tilde{g}:\Lambda_{\tilde{g}}\mapsto\Lambda_{\tilde{g}}$ is $C^1$-conjugated to a diffeomorphism $C^1$ close to the dynamics  $f:\mathbb{T}^3\mapsto\mathbb{T}^3$. Then the bundle $E^c|_{\Lambda_{\tilde{g}}}$ is volume expanding and there is a periodic orbit of index one with complex eigenvalues along the bundle $E^c$ contained in $\Lambda_{\tilde{g}}$. Hence the bundle $E^c|_{\Lambda_{\tilde{g}}}$ also has no finer dominated splitting. As a consequence,  we have the following conclusion.

 \begin{lemma}
 For generic diffeomorphism in $\mathcal{V}$, there is a partially hyperbolic homoclinic class, such that any ergodic measure supported on it has at least one positive Lyapunov exponent along $E^c$.
 \end{lemma}

\vskip 5mm
{\bf Acknowlegement: }
We would like to  thank Christian Bonatti and Sylvain Crovisier  whose comments and suggestions help us to get stronger result than the original version of this paper.

We would also thank Lan Wen, Shaobo Gan and Dawei Yang for useful comments.

Jinhua Zhang  would like to thank China Scholarship Council (CSC) for financial support (201406010010).

%%%%%%%%%%%%%%%%%%%%%%%%%%%%%%%%%%%%%%%%%%%%%%%%%%%%%%%%%%%%%%%%%%%%%%
\bibliographystyle{plain}

\begin{thebibliography}{99}


%\bibitem[ABC]{ABC} F. Abdenur, Ch. Bonatti and S. Crovisier,
%\emph{Nonumiform hyperbolicity for $C^{1}$ generic diffeomorphisms}. Israel J. Math. 183 (2011), 1--60.


\bibitem[ABCDW]{ABCDW} F. Abdenur, Ch. Bonatti, S. Crovisier, L. J. D\'iaz and L. Wen,  \emph{Periodic points and homoclinic classes}. Ergodic Theory $\&$ Dynamical System 27 (2007), 1--22.

\bibitem [AS]{AS} R. Abraham and S. Smale, \emph{Nongenericity of $\Omega$-stablity}. Global Analysis I, Proc. Symp. Pure Math. 1968, AMS 14 (1970), 5--8.

\bibitem[BB]{BB} J. Bochi and Ch. Bonatti, \emph{Perturbation of the Lyapunov spectra of periodic orbits}. Proc. Lond. Math. Soc. 105 (2012), 1--48.


\bibitem[BBD1]{BBD1} J. Bochi, Ch. Bonatti and L. J. D\'iaz, \emph{Robust vanishing of all Lyapunov exponents for iterated function systems}. Math. Z. 276 (2014), 469--503.

\bibitem[BBD2]{BBD2} J. Bochi, Ch. Bonatti and L. J. D\'iaz, \emph{Robust criterion for the existence of nonhyperbolic ergodic measures}.	 ArXiv:1502.06535.

\bibitem[B]{B} Ch. Bonatti, \emph{A semi-local argument for robust transitivity}.  Lecture at IMPA's Dynamical
Systems Seminar, August, 1996.


\bibitem[BC]{BC} Ch. Bonatti and S. Crovisier, \emph{R$\acute{e}$currence et g$\acute{e}$n$\acute{e}$ricit$\acute{e}$}. Invent. Math. 158 (2004), 33--104.

\bibitem[BCDG]{bcdg} Ch. Bonatti, S. Crovisier, L. J. D\'iaz and N. Gourmelon, \emph{Internal perturbations of homoclinic classes: non-domination, cycles, and self-replication.} Ergodic Theory $\&$ Dynamical Systems 33 (2013), 739--776.

\bibitem[BDG]{BDG} Ch. Bonatti, L. J. D\'iaz and A. Gorodetski, \emph{Non-hyperbolic ergodic measures with large support}. Nonlinearity 23
    (2010), 687--705.

\bibitem[BDP]{BDP} Ch. Bonatti, L. J. D\'iaz and E. R. Pujals, \emph{A $C^{1}$ generic dichotomy for diffeomorphisms: Weak forms of
hyperbolicity or infinitely many sinks or sources}. Ann. of Math. 158 (2003), 355--418.

\bibitem[BDPR]{BDPR}Ch. Bonatti, L. J. D\'iaz, E. R. Pujals and J. Rocha, \emph{Robustly transitive sets and heterodimensional cycles}. Ast$\acute{e}$risque 286 (2003), 187--222.

\bibitem[BDV]{BDV} C. Bonatti, L. D\'iaz and M. Viana, \textit{Dynamics beyond uniform hyperbolicity. A global geometric and probabilistic perspective.} Encyclopaedia of Mathematical Sciences 102. Mathematical Physics, III. Springer-Verlag (2005).

\bibitem[BV]{BV} Ch. Bonatti and M. Viana, \emph{SRB measures for partially hyperbolic systems whose central direction is mostly
contracting}. Israel J. Math. 115 (2000), 157--193.

\bibitem[CLR]{CLR}Y. Cao, S. Luzatto and I. Rios, \emph{Some non-hyperbolic systems with strictly non-zero Lyapunov exponents for all invariant measures: horseshoes with internal tangencies}. Disc. Cont. Dyn. Syst. 15 (2006), 61--71.


\bibitem[CCGWY]{CCGWY} Ch. Cheng, S. Crovisier, S. Gan, X. Wang and D. Yang, \emph{Hyperbolicity versus non-hyperbolic ergodic measures inside homoclinic classes}.	 ArXiv:1507.08253.

%\bibitem[Co]{conley} C. Conley, \emph{Isolated Invariant Sets and the Morse Index.} CBMS Regional Conference Series in Mathematics 38,
%American Mathematical Society, Providence. RI (1978).

% \bibitem[C]{C} S. Crovisier, \emph{Partial hyperbolicity far from homoclinic bifurcations}. Advance in Math. 226 (2011), 673--726.

\bibitem[DG]{DG} L. J. D\'iaz and A. Gorodetski, \emph{Non-hyperbolic ergodic measures for non-hyperbolic homoclinic classes}. Ergodic Theory $\&$ Dynamical Systems 29 (2009), 1479--1513.

\bibitem[G]{G} S. Gan, \emph{A necessary and sufficient condition for the existence of dominated splitting with a given index}. Trends in Mathematics 7 (2004), 143--168.

\bibitem[Go]{Go} N. Gourmelon, \emph{A Frank's lemma that preserves invariant manifolds}. To appear at Ergodic Theory $\&$ Dynamical Systems. ArXiv:0912.1121.

 \bibitem [GIKN]{GIKN} A. Gorodetski, Yu. S. Ilyashenko, V. A. Kleptsyn and M. B. Nalsky,  \emph{Nonremovability of zero Lyapunov exponents}. (Russian) Funktsional. Anal. i Prilozhen. 39 (2005), 27--38; translation in Funct. Anal. Appl. 39 (2005), 21--30.

%\bibitem[H]{H}S.Hayashi, \emph{Connecting invariant manifolds and the solution of $C^1$-stability and $\Omega$-stability conjectures
% for flows}. Ann. of Math, 145 (1997), 81--137.

\bibitem[HPS]{HPS} M. Hirsch, C. Pugh and M. Shub, \emph{Invariant manifolds}. Lecture Notes in Math. 583. (1977).

\bibitem [K]{K} A. Katok,  \emph{Lyapunov exponents, entropy and periodic orbits for diffeomorphisms}.
Publ. Math. I.H.E.S, 51 (1980), 137--173.

\bibitem[KN]{KN} V. A. Kleptsyn and M. B. Nalsky,  \emph{Stability of the existence of nonhyperbolic measures for $C^1$ diffeomorphisms}.
 Funktsional. Anal. i Prilozhen. 41 (2007), 30--45; translation in Funct. Anal. Appl.41(2007), 271--283.

\bibitem[O]{O} V. I. Oseledets, \emph{A multiplicative ergodic theorem: Lyapunov characteristic numbers for dynamical systems}. Trans. Moscow. Math. Soc. 19 (1968), 197--231.

\bibitem[P]{pesin} Y. Pesin, \emph{Characteristic Lyapunov exponents and smooth ergodic theory.} Usp. Mat. Nauk. 32 (1977), 55--112.


%\bibitem[Pl]{p} V. Pliss, \emph{On a conjecture of Smale.} Differ. Uravn. 8 (1972), 262--268.

%\bibitem[PR]{PR}C. Pugh and C. Robinson, \emph{The $C^1$ closing lemma, including Hamiltonians}. Ergodic Theory $\&$ Dynamical Systems 3 (1983), 261--313.

%\bibitem[PS]{PS} E. Pujals and M. Sambarino, \emph{Homoclinic tangencies and hyperbolicity for suface diffeomorphisms.} Ann. of Math. 151 (2000), 961--1023.

%\bibitem[W]{W} P. Walter, \emph{An introduction to ergodic theory}. Graduate Texts in Mathematics 79, Springer, New York-Berlin, 1982.

\bibitem[Wa]{wang} X. Wang, \emph{Hyperbolicity versus weak periodic orbits inside homoclinic classes.}  ArXiv:1504.03153.

%\bibitem[WX]{WX} L. Wen and Z. Xia, \emph{$C^1$ connecting lemmas}. Trans. Amer. Math. Soc. 352, 11(2000), 5213--5230.

\end{thebibliography}
%%%%%%%%%%%%%%%%%%%%%%%%%%%%%%%%%%%%%%%%%%%%%%%%%%%%%%%%%%%%%%%%%%%%%%%%%

%%%%%%%%%%%%%%%%%%%%%%%%%%%%%%%%%%%%%%%%%%%%%%%%%%%%%%%%%%%%%%%%%%%%%%%%%%%%%%
\vskip 2mm

\noindent Xiaodong Wang,

\noindent{\small School of Mathematical Sciences\\
Peking University, Beijing 100871, China}\\
\noindent and\\
\noindent{Laboratoire de Math\'ematiques d'Orsay\\
 Universit\'e Paris-Sud 11, Orsay 91405, France}

\noindent {\footnotesize{E-mail : xdwang1987@gmail.com }}

\vskip 2mm
\noindent Jinhua Zhang,

\noindent{\small School of Mathematical Sciences\\
Peking University, Beijing 100871, China}

\noindent {\footnotesize{E-mail :zjh200889@gmail.com }}\\
\noindent and\\
\noindent {\small Institut de Math\'ematiques de Bourgogne\\
UMR 5584 du CNRS}

\noindent {\small Universit\'e de Bourgogne, 21004 Dijon, France}\\
\noindent {\footnotesize{E-mail : jinhua.zhang@u-bourgogne.fr}}

\end{document}